\documentclass{article}

\usepackage{amsmath,amsthm,amssymb}

\usepackage{stackrel}

\newtheorem{theorem}{Theorem}[section]

\newtheorem{corollary}[theorem]{Corollary}

\newcommand{\NN}{\mathbb{N}}
\newcommand{\ZZ}{\mathbb{Z}}

\newcommand{\BB}{\mathbb{B}}

\newcommand{\fa}[1]{\forall #1 \,}
\newcommand{\ph}{\varphi}
\newcommand{\limplies}{\rightarrow}

\title{Character and object\footnote{This is an expanded version of the corresponding journal publication. In the latter, Section~\ref{metaphysics:section} has been shortened, and Sections~\ref{frege:section} and \ref{frege:section:b} are omitted. This work draws on the second author's Carnegie Mellon MS thesis \cite{morris:11}. We are grateful to Michael Detlefsen and the participants in his \emph{Ideals of Proof} workshop, which provided helpful feedback on portions of this material in July, 2011. We are also grateful to Paddy Blanchette, Roy Cook, Jeremy Heis, Pen Maddy, Marco Panza, and Marcus Rossberg. Avigad's work has been partially supported by National Science Foundation grant DMS-1068829 and Air Force Office of Scientific Research grant FA9550-12-1-0370.}}

\author{Jeremy Avigad and Rebecca Morris}

\begin{document}

\maketitle

\begin{abstract}
 In 1837, Dirichlet proved that there are infinitely many primes in any arithmetic progression in which the terms do not all share a common factor. Modern presentations of the proof are explicitly higher-order, in that they involve quantifying over and summing over \emph{Dirichlet characters}, which are certain types of functions. The notion of a character is only implicit in Dirichlet's original proof, and the subsequent history shows a very gradual transition to the modern mode of presentation. 
 
 In this essay, we describe an approach to the philosophy of mathematics in which it is an important task to understand the roles of our ontological posits and assess the extent to which they enable us to achieve our mathematical goals. We use the history of Dirichlet's theorem to understand some of the reasons that functions are treated as ordinary objects in contemporary mathematics, as well as some of the reasons one might want to resist such treatment. We also use these considerations to illuminate the formal treatment of functions and objects in Frege's logical foundation, and we argue that his philosophical and logical decisions were influenced by many of the same factors.
\end{abstract}

\tableofcontents

\section{Introduction}
\label{introduction:section}

The philosophy of mathematics has long been concerned with the nature of mathematical objects, and the proper methods for acquiring mathematical knowledge. But as of late some philosophers of mathematics have begun to raise questions of a broader epistemological character: What does it mean to properly \emph{understand} a piece of mathematics? In what sense can a proof be said to \emph{explain} a mathematical fact? In what senses can one proof be viewed as better than another one that establishes the same theorem? What makes a concept fruitful, and what makes one definition more natural than another? Why are certain historical developments viewed as important advances? Questions like these are sometimes classified as pertaining to the \emph{methodology} of mathematics, in contrast to more traditional ontological concerns.

One of our goals in this essay is to argue that methodology and ontology cannot be so cleanly separated. Certainly part of the justification for our ontological commitments stems from the positive effects those commitments have on the practice, and, conversely, ``internal'' methodological shifts are influenced by a broader conception as to what is permissible. In Section~\ref{metaphysics:section}, we describe a model for historical change that closely links ontological and methodological considerations.

One of the hallmarks of the nineteenth century transition to modern mathematics was the adoption of implicit or explicit set-theoretic language and methods. For Gauss \cite{gauss:01}, the number-theoretic relation of congruence modulo $m$ was a relation that was similar to equality, and addition and multiplication modulo $m$ were operations on integers that respect that relation. Today, however, we can form the quotient structure of integers modulo $m$, which consists of classes of integers that are equivalent modulo $m$. Addition and multiplication then lift to operations on these classes. This amounts to \emph{reifying} the property of being equivalent to an integer $a$ modulo $m$ to an object, $[a]$, the equivalence class of $a$. Similarly, to restore the property of unique factorization to the algebraic integers in a cyclotomic field, Kummer \cite{kummer:46} introduced properties $P(\alpha)$ that were meant to be interpreted as the assertion that $\alpha$ is divisible by a certain ``ideal divisor.'' Dedekind \cite{dirichlet:63b} later reified the property $P$ to the class of $\alpha$ that satisfy it, thereby giving rise to the modern notion of an ideal in a ring of integers. Other nineteenth century examples include the construction of quotient groups, or the lifting of Gauss' operation of ``composition'' of binary quadratic forms to equivalence classes of such forms.

What these instances have in common is that they involve treating certain higher-order entities --- classes of integers, classes of algebraic integers, classes of quadratic forms, or classes of elements in a group or a ring --- as objects in their own right. By this we mean that, in particular, one can quantify over them, sum over them, and define operations on them. Moreover, one can consider algebraic structures whose elements are such classes, much as one can consider algebraic structures whose elements are integers or real or complex numbers.

Much of what can be said about the treatment of classes as objects in the nineteenth century applies to the treatment of functions as objects as well. In 1837, Dirichlet proved that there are infinitely many prime numbers in any arithmetic progression in which the terms do not all share a common factor. Our goal here is to study the role that certain types of functions, called \emph{Dirichlet characters}, play in contemporary presentations of Dirichlet's proof, and the historical process that has led to our contemporary understanding.

In Section~\ref{metaphysics:section}, we present a framework for assessing the ontological commitments of a body of mathematics, one which is informed by, and can inform, the history of mathematics. In Section~\ref{overview:section}, we provide an overview of Dirichlet's proof, and in Section~\ref{functions:section}, we clarify the senses in which contemporary presentations treat characters as ordinary mathematical objects. Despite the name, the notion of a Dirichlet character is not present in Dirichlet's original presentation. In Sections~\ref{dirichlet:section} and~\ref{transition:section}, we describe the history of presentations of Dirichlet's theorem, which shows a fitful and gradual transition to modern terminology and usage. In doing so, we draw on a detailed historical study that we have carried out in another work \cite{avigad:morris:unp}, which we will refer to as ``Concept'' in the presentation below.

In Section~\ref{analysis:section}, we argue that, as per the model presented in Section~\ref{metaphysics:section}, the gradual adoption of the modern treatment of characters is best viewed as an ontological response to pragmatic mathematical concerns, and we explore some of the considerations that bear on the rationality of the outcome. Thus we use the history to help us understand and assess some of the reasons that we treat functions as objects in current mathematical practice. Complementing the mathematical narrative, in Section~\ref{frege:section}, we consider Frege's conflicted attitudes towards the treatment of functions as objects, and in Section~\ref{frege:section:b}, we argue that key choices in the design of his formal system were motivated by the same sorts of considerations. This is not to say that Frege's logico-philosophical concerns should be seen as properly mathematical, or vice versa. Rather, they both stem from the need to balance two key desiderata: the desire, on the one hand, for flexible and uniform ways of dealing with higher-order entities in the many guises in which they appear, and the desire, on the other hand, to make sure that the methods of doing so are clear, coherent, and meaningful.

\section{From methodology to ontology}
\label{metaphysics:section}

Let us start by distinguishing between two kinds of questions one can ask, having to do with the existence of mathematical objects. On the one hand, we can ask question such as:
\begin{itemize}
\item Is there a nontrivial zero of the Riemann zeta function whose real part is not equal to $1/2$?
\item Are there noncyclic simple groups of odd order?
\end{itemize}
These are fundamentally \emph{mathematical} questions. Answering them is not easy: the Riemann hypothesis posits a negative answer to the first, while the Feit-Thompson theorem, a landmark in finite group theory, provides a negative answer to the second. But even in the first case, where we do not know the answer to the question, we feel that we have a clear sense as to what kind of argument would settle the issue one way or another. Put simply, questions like these can be addressed using conventional mathematical methods. In contrast, there are questions like these:
\begin{itemize}
 \item Do the natural numbers (really) exist, and what sorts of things are they?
 \item Are there infinite totalities?
 \item What kinds of sets and functions exist (if any), and what properties do they have?
 \item Are there infinitesimals, fluxions, fluents, and ultimate ratios?
\end{itemize}
These are questions as to the ultimate nature of mathematics and its objects of study, and seem to call for a more general, open-ended \emph{philosophical} analysis. What is sought is not just an axiomatization of mathematics or an enumeration of the mathematical objects that exist, but also explanation as to why we are justified in asserting in their existence, with an overall account that squares with broader epistemological and scientific concerns.

The distinction between the two types of questions may call to mind the logical positivists' distinction between questions that are ``internal'' to a linguistic framework, and ``external'' or ``pragmatic'' questions pertaining to the choice of a framework itself. Some take this distinction to be have been repudiated, decisively, by the criticisms of W.~V.~O.~Quine \cite{quine:51}. But keep in mind that Quine's arguments, which were directed against the claim that there is a sharp, principled distinction between the two sorts of questions, were not meant to show that there is no difference between them at all. In locating both kinds of questions on the common continuum of scientific inquiry, he did not deny that different kinds of questions call for different sorts of answers; indeed, his influential \emph{Word and Object} \cite{quine:60} is an extended exploration of the considerations that he took to bear on ``philosophical'' questions of the latter sort. Nothing we say below commits us to a sharp distinction, and it seems relatively uncontroversial to say that insofar as any rational arguments can be brought to bear on the second group of questions, these will look different from the kinds of arguments that are brought to bear on the first.

Despite their different characterizations of the philosophical project, Carnap and Quine shared the view that ontological questions come down to pragmatic questions as to the choice of a conceptual framework. Here is what Carnap had to say about our scientific commitments to abstract objects:
\begin{quote}
 The acceptance cannot be judged as being either true or false because it is not an assertion. It can only be judged as being more or less expedient, fruitful, conducive to the aim for which the language is intended. Judgments of this kind supply the motivation for the decision of accepting or rejecting the kind of entities.\footnote{Here and below, when the bibliographic entry of a work includes a reprinted version, page numbers in the references refer to the reprinted version. Similarly, when the bibliographic entry includes an English translation, our translations are taken from that source, unless we indicate otherwise. Where no translation is listed, the translations are our own. The original versions of most of the mathematical sources quoted here can be found in ``Concept,'' so we have not reproduced them here.} \cite[p.~250]{carnap:50}
\end{quote}
Quine offered the following amendment:
\begin{quote}
 Consider the question whether to countenance classes as entities. This, as I have argued elsewhere, is the question whether to quantify with respect to variables which take classes as values. Now Carnap has maintained that this is a question not of matters of fact but of choosing a convenient language form, a convenient conceptual scheme or framework for science. With this I agree, but only on the proviso that the same be conceded regarding scientific hypotheses generally. \cite[p.~43]{quine:51}
\end{quote}

We take these views seriously here, seeing it as an important philosophical task to clarify the role of our ontological posits with respect to ordinary mathematical activity, and evaluate their efficacy towards achieving our mathematical goals. This amounts to something like the naturalist approaches to the philosophy of mathematics advocated by Kitcher \cite{kitcher:88}, Burgess \cite{burgess:08}, and Maddy \cite{maddy:97}, focused on specific aspects of mathematical practice.

% On this view, philosophy becomes a kind of means-ends analysis. It is therefore important to keep in mind that scientific ends are usually best described relative to a scientific context. While a mathematician may have good reason to admit certain abstract objects with all their attendant properties, a logician may have good reason to adopt a more minimal ontology when it comes to studying the metamathematical effects of our axiomatic assumptions; and there may be good reasons for psychologists to avoid reference to abstract objects altogether in an account of how children learn to grasp mathematical concepts. It seems fruitless to grant any of these contexts the final say as to whether mathematical objects ``really'' exist; the objects in question play distinct roles in distinct theories with distinct explanatory and predictive goals. In situations where the scientific contexts overlap, some translation and reinterpretation is inevitable, and part of the metaphysician's task is to help keep these cross-disciplinary collaborations running smoothly. Our view does not rule out hope for an overarching frame that is broad enough to fit together the various scientific snapshots; but we do insist that, ultimately, the success of such a metaphysics needs to be judged in terms of its contribution to the pragmatic needs of the working scientist.

How, then, should such an analysis proceed? It is instructive to consider those historical situations in which the mathematical community faced possibilities for methodological or ontological expansion and reacted accordingly. For example, it is helpful to consider the ancient Greek idealizations of number and magnitude, and the theory of proportion; the gradual acceptance of negative numbers, and then complex numbers, in the Western tradition; the use of algebraic methods in geometry, infinitesimals in the calculus, points at infinity in projective geometry; the development of the function concept from Euler to modern times; the gradual set-theoretic treatment of algebraic objects like cosets, ideals, equivalence classes in the nineteenth century; and so on. By studying the historical concerns regarding these expansions as well as the pressures that led to their ultimate acceptance, we can hope to better understand the factors that influence such developments.

Indeed, at junctures like these, historical developments tend to follow a common pattern. First, expansions are met with resistance, or at least, extreme caution. Sometimes, the expansions can be explained in terms of the more conservative practice; for example, complex numbers can be interpreted as ordered pairs, algebraic solutions to geometric problems can be reinterpreted geometrically, and equations can be rewritten to avoid consideration of negative quantities. In other cases, the expansions are not generally conservative, but, at least, can be explained away in particular instances; for example, arguments involving infinitesimals can sometimes be interpreted in terms of ``ultimate ratios'' in a geometric diagram, and operations on abstract objects can sometimes be understood as operations on explicit representations. This makes it possible to adopt the expansions, tentatively, as convenient shorthand for more tedious but conservative arguments. Over time, the rules and norms that govern the expansions are clarified, and the expansions themselves prove to be convenient, or even indispensable, while they do not cause serious problems. Over time, the mathematical community grows used to them, to the point where they become part of the usual business of mathematics.

Whiggish narratives tend to dismiss such historical hand-wringing and shilly-shallying as short-sighted conservativism that stands in the way of mathematical progress. We, however, prefer to view it as a rational response to the proposed expansions, whereby the benefits are carefully weighed against the concerns. In hindsight, we tend to make too little of the pitfalls associated with an ontological or methodological expansion. To start with, there are concerns about the \emph{consistency} and \emph{coherence} of the new methods, that is, worries as to whether the changes will lead to mistakes, false results, or utter nonsense, perhaps when employed in situations that have not even been imagined. Kenneth Manders has also emphasized the importance of maintaining \emph{control} of our mathematical practices \cite{manders:08}. Mathematics requires us to be able to come to agreement as to whether a proof is correct, or whether a given inference is valid or not. If new objects come with rules of use that are not fully specified, or vague, or unclear, the practice is in danger of breaking down. In a sense, this concern comes prior to concerns of consistency: if it is not clear what properties abstract magnitudes, negative numbers, complex numbers, infinitesimals, sets, and ``arbitrary'' functions have, it doesn't even make sense to ask whether using them correctly will lead to contradictions.\footnote{Mathematics, however, often gets by surprisingly well with concepts that are problematic, incompletely specified, and not fully understood, something which has been emphasized by Wilson \cite{wilson:94} and Urquhart \cite{urquhart:08}.} 

And then there are further concerns as to whether the new methods are \emph{meaningful} and \emph{appropriate} to mathematics. Even if a body of methods is consistent and clearly specified, it may still fail to provide us with the results we are after. If we expect an existence proof to yield certain kinds of information about the object that is asserted to exist, methods that fail to provide that sort of information do not constitute mathematics---or, at least, not the kind of mathematics we should be doing. If you expect a mathematical theory to make scientific predictions that we can act on rationally, it is a serious concern as to whether the new methods can deliver.

In short, the concerns are not easily set aside. What, then, are the factors that might sway a decision in favor of an expansion? Mathematicians tend to wax poetic in their praise of conceptual advances, highlighting the power of new methods, the elegance and naturality of the resulting theory, and the insight and depth of the associated ideas. Part of our goal here is to de-romanticize these virtues and gain clarity as to what might be achieved. In many instances, the virtues in question have a lot to do with efficiency and economy of thought:\footnote{The phrase is borrowed from Ernst Mach's \emph{The Science of Mechanics} \cite{mach:93}; we are grateful to Michael Detlefsen for bringing this to our attention.} we tend to value methods that make it possible to solve problems that were previously unsolvable, or simplify proofs and calculations that were previously tedious, complex, and error-prone. Below we will consider specific ways in which ontological and methodological expansions help us manage complex tasks by suppressing irrelevant detail, making key features of a problem salient, and keeping key information ready-to-hand. We will also try to understand the way they make it possible to generalize and extend results, and facilitate the transfer of ideas to other domains.

To summarize our high-level historical model: when mathematics is faced with methodological expansion, benefits such as simplicity, generality, and efficiency are invariably weighed against concerns as to the consistency, cogency, and appropriateness of the new methods. Sufficient benefit encourages us to entertain the changes cautiously, while trying to minimize the dangers involved. Cogency is obtained by working out the norms and conventions that govern the new methods. Consistency may not be guaranteed, but our experiences over time can bolster our faith that the new methods do not cause problems. In this regard, initial checks that the new methods are partially conservative over the old ones helps preserve mathematical meaning, and reassures us that even if the new methods turn out to be problematic, one will be able to restrict their scope in such a way that preserves their utility.\footnote{Wittgenstein's discussion of contradiction is interesting in this regard; see \cite[Lectures XI--XII]{wittgenstein:89}.} The philosophy of mathematics should give us better means to evaluate such expansions: to talk about the cogency of a mathematical argument and whether it delivers the desired result, and to understand the ways in which our ontological posits and methodological expansions improve our ability to reason effectively.

A salient feature of our approach is that we aim to take mathematics at face value: when our best mathematical theories tell us that numbers and functions exist, our best philosophical theories should not repudiate those claims. This feature is common to other approaches to the philosophy of mathematics, such as the platonism espoused by William Tait \cite[\S 5]{tait:86}, or the non-eliminative structuralism proposed by Charles Parsons \cite[\S 18]{parsons:08}. On the other hand, we recognize an ongoing need for sustained reflection on our mathematical goals and methods, in order to better understand and improve that practice. Taking mathematics at face value doesn't mean viewing it as fixed and unchanging; mathematics has evolved for centuries and will continue to do so, guided, we hope, by thoughtful reflection of this sort.

At times, it may seem that our treatment of ontological questions verges on a kind of formalist instrumentalism, for example, the view that there is nothing more to mathematics than linguistic conventions, which are to be adjudicated on the basis of ``pragmatic'' concerns. To be sure, we take pragmatic concerns to be an important target of philosophical study, but insofar as there is anything to be made of the realism/anti-realism debate with respect to mathematics, nothing we say here should preclude a realist position. For example, Hilary Putnam has argued that
\begin{quote}
\ldots at least when it comes to the theories that scientists regard as most fundamental\ldots we should regard all of the rival theories as candidates for truth or approximate truth, and that \emph{any philosophy of mathematics that would be inconsistent with so regarding them should be rejected}. \cite[p.~184]{putnam:12}
\end{quote}
Moreover:
\begin{quote}
\ldots a \emph{prima facie} attractive position---realism with respect to the theoretical entities postulated by physics, combined with \emph{antirealism} with respect to mathematical entities and/or modalities---doesn't work. [\emph{ibid.}, p.~188]
\end{quote}
Considering our mathematical and scientific theories as ``candidates for truth or approximate truth'' does not preclude reflecting on those theories and bringing pragmatic considerations to bear on the choices among them. Indeed, that is an integral part of the scientific enterprise, and it is the kind of activity we hope to support.

\section{An overview of Dirichlet's theorem}
\label{overview:section}

Two integers, $m$ and $k$, are said to be \emph{relatively prime}, or \emph{coprime}, if they have no common factor. In 1837, Dirichlet proved the following:
\begin{theorem}
\label{dirichlet:theorem}
 If $m$ and $k$ are relatively prime, the arithmetic progression $m, m + k, m + 2k, \ldots$ contains infinitely many primes.
\end{theorem}
In other words, if $m$ and $k$ are relatively prime, there are infinitely many primes congruent to $m$ modulo $k$.

In 1798, Legendre had assumed this, without justification, in a purported proof of the law of quadratic reciprocity. Gauss pointed out this gap, and presented two proofs of quadratic reciprocity in his \emph{Disquisitiones Arithmeticae} of 1801, which do not rely on that fact. He ultimately published six proofs of quadratic reciprocity, and left two more in his \emph{Nachlass}, but he never proved the theorem on primes in an arithmetic progression. Dirichlet's proof is notable not only for settling a longstanding open problem, but also for its sophisticated use of analytic methods to prove a number-theoretic statement.

\subsection{Euler's proof that there are infinitely many primes}

As Dirichlet himself made clear, the conceptual starting point for his proof lies in the work of Euler. In the \emph{Elements}, Euclid proved that there are infinitely many primes, but his proof does not provide much information about how they are distributed. Euler, in his \emph{Introductio in Analysin Infinitorum} \cite{Euler48}, proved the following:
\begin{theorem}
\label{euler:thm}
The series $\sum_{q}\frac{1}{q}$ diverges, where the sum is over all primes $q$.
\end{theorem}
This implies that there are infinitely many primes, but also says something more about their density. For example, since we know that the series $\sum_{n}\frac{1}{n^{2}}$ is convergent, it tells us that, in a sense, there are ``more'' primes than there are squares.

Euler's proof of Theorem~\ref{euler:thm} centers around his famous zeta function,
\[
\zeta(s) = \sum_{n=1}^{\infty}n^{-s} = 1 + \frac{1}{2^s} + \frac{1}{3^s} + \ldots,
\]
defined for a real variable $s$. (The zeta function was later extended by Riemann to the entire complex plane via analytic continuation.) It is not hard to show that the series $\zeta(s)$ converges whenever $s > 1$. In that case, the infinite sum can also be expressed as an infinite product:
\begin{equation}
\label{euler:product:eqn}
 \sum_{n=1}^{\infty} n^{-s}= \prod_{q}\left(1-\frac{1}{q^{s}}\right)^{-1},
\end{equation}
where the product is over all primes $q$. This is known as the \emph{Euler product formula}. Roughly, this holds because we can write each term of the product as the sum of a geometric series,
\[
\left(1-\frac{1}{q^{s}}\right)^{-1} = 1 + q^{-s} + q^{-2s} + \ldots
\]
and then expand the product into a sum. The unique factorization theorem tells us that every integer $n > 1$ can be written uniquely as a product $q_1^{i_1} \cdot q_2^{i_2} \cdots q_k^{i_k}$. This means that the term $n^{-s} = q_1^{-i_1 s} \cdot q_2^{-i_2 s} \cdots q_k^{-i_k s}$ will occur exactly once in the expansion, corresponding to the choice of the $i_j$th element of the sum for each $q_j$, and the choice of $1$ in every other sum. Since we are dealing with infinite sums and products, the Euler product formula implicitly makes a statement about limits, and some care is necessary to make the argument precise; but this is not hard to do.

If we take the logarithm of each side of the product formula and appeal to properties of the logarithm function, we obtain
\[
 \log\sum_{n=1}^{\infty}n^{-s}=\sum_{q}-\log\left(1-\frac{1}{q^{s}}\right).
\]
Using the Taylor series expansion
\[
\log(1 - x) = -x - x^2 / 2 - x^3 / 3 - \ldots
\]
and changing the order of summations yields
\[
\log\sum_{n=1}^{\infty}n^{-s}=\sum_{q}\frac{1}{q^{s}} + \sum_{n=2}^{\infty}\frac{1}{n}\sum_{q}\frac{1}{q^{ns}}.
\]
Remember that we want to show that $\sum_{q}\frac{1}{q}$ diverges. Notice that the first term on the right-hand side of the above equation is $\sum_{q}\frac{1}{q^{s}}$. Thus we should consider what happens as $s$ tends to 1 from above. One can show that the second term on the right-hand side is bounded by a constant that is independent of $s$, a fact that can be expressed using ``big O'' notation as follows:
\begin{equation}
\label{euler:primes:eqn}
\log\sum_{n=1}^{\infty}\frac{1}{n^s}=\sum_{q}\frac{1}{q^{s}} + O(1).
\end{equation}
As $s$ approaches $1$ from above, the left-hand side clearly tends to infinity. Thus, the right-hand side, $\sum_{q}\frac{1}{q^{s}}$, must also tend to infinity, which implies that $\sum_{q}\frac{1}{q}$ diverges.

\subsection{Dirichlet's approach}
\label{dirichlet:approach:section}

To make the ideas more perspicuous, Dirichlet first considered Theorem~\ref{dirichlet:theorem} in the special case where the common difference is a prime number $p$. Any prime $q$ other than $p$ leaves a remainder of $1, \ldots, p-1$ when divided by $p$. Splitting up the sum in (\ref{euler:primes:eqn}) we then have
\begin{equation}
\label{euler:primes:eqn:b}
\log\sum_{n=1}^{\infty}\frac{1}{n^s}=\sum_{q \equiv 1 \bmod p} \frac{1}{q^{s}} + \sum_{q \equiv 2 \bmod p} \frac{1}{q^{s}} + \ldots + \sum_{q \equiv p-1 \bmod p} \frac{1}{q^{s}} + O(1).
\end{equation}
This shows that (\ref{euler:primes:eqn}) is too crude to prove Theorem~\ref{dirichlet:theorem}: to show that there are infinitely many primes congruent to $m$ modulo $p$, we need to show that the $m$th term on the right-hand side tends to infinity, not just the sum of all such terms. More work is therefore needed to tease apart the contribution of the primes modulo $m$, for each nonzero residue $m$ modulo $p$.

Dirichlet sketched his proof in a three-page note announcing the result in 1837 \cite{Dirichletshort}, before spelling out the details in a later publication \cite{Dirichlet37}. The method relies on a trick that seems to come out of nowhere. We describe the trick here, and in the Appendix offer an explanation as to how Dirichlet may have come upon this approach. 

It is a fact from number theory that for any prime number $p$, there is a number $g$, such that the powers $g^0, g^1, g^2, \ldots, g^{p-2}$ modulo $p$ are exactly the nonzero residues $1, 2, 3, \ldots, p-1$ modulo $p$ in some order. Such an element $g$ is called a \emph{primitive root modulo $p$}. For example, when $p = 11$, we can choose $g = 2$. In that case, the powers of $g$ modulo 11 are
\[
 1, 2, 4, 8, 5, 10, 9, 7, 3, 6, 
\]
which are just the numbers from 1 to 10 listed in a different order. Notice that the next element on the list would be $1$ again, and the list cycles. In general, if $g$ is a primitive root modulo $p$, then $g^{p-1}$ is equal to $1$ modulo $p$.

The statement that $g$ is a primitive root modulo $p$ means that for each nonzero residue $m$ modulo $p$, there is an exponent $\gamma$ between $0$ and $p-2$, with the property that $g^\gamma$ is equal to $m$ modulo $p$. We will denote this exponent $\gamma_m$ and call it the \emph{index} of $m$ modulo $p$ with respect to $g$, as Dirichlet did. For example, consulting the list above, we see that the index of $10$ is $5$, because $2^5$ is equal to $10$ modulo $11$. The function $n \mapsto \gamma_n$ behaves like a logarithm, in the sense that if $m$ and $n$ are nonzero residues modulo $p$, $\gamma_{mn}$ is equal to $\gamma_m + \gamma_n$ modulo $p - 1$. This is because we have
\[
 g^{\gamma_m + \gamma_n} = g^{\gamma_m} g^{\gamma_n} = m n \bmod p,
\]
and so $\gamma_m + \gamma_n$ modulo $p - 1$ is the exponent corresponding to $mn$. 

We now turn our attention from integer roots modulo a prime to the notion of a complex root of unity. In general, if $n$ is any integer, the equation $x^n = 1$ will have $n$ distinct roots in the complex numbers. Moreover, we can choose such a root, $\omega$, that is primitive in the sense that $\omega^0, \omega^1, \omega^2, \ldots, \omega^{n-1}$ are all such roots; taking $\omega = e^{2 \pi i / n}$ will do. 

Notice that we are now using the phrase ``primitive root'' in two distinct, but related, senses: to refer to primitive roots modulo a prime, and to refer to primitive roots of unity. For future reference, notice also that the expression $x^n - 1$ factors as $(x - 1) (x^{n - 1} + \ldots + x^2 + x+ 1)$. So, for any complex number $x$, if $x$ is a solution to $x^n = 1$ other than $1$, we have $x^{n - 1} + \ldots + x^2 + x+ 1 = 0$.

Returning to Dirichlet's theorem, let $p$ be any prime, fix a primitive root $g$ modulo $p$, and let $\omega$ be any $(p - 1)$st root of 1, primitive or not. Consider the function $\chi(n)$ which maps any nonzero residue $n$ to the value $\omega^{\gamma_n}$. The function $\chi$ is \emph{multiplicative}, which is to say, $\chi(mn) = \chi(m) \chi(n)$ for any two nonzero residues $m$ and $n$. This holds because
\[  
\chi(mn) = \omega^{\gamma_{mn}} = \omega^{\gamma_m + \gamma_n} = \omega^{\gamma_n} \omega^{\gamma_n} = \chi(m) \chi(n).
\]
In the next section, we will see the functions $\chi$ are exactly the \emph{characters} on the group of nonzero residues modulo $p$. Here, following Dirichlet, we will avoid writing $\chi(n)$ and stick with the notation $\omega^{\gamma_n}$.

A crucial ingredient in Dirichlet's proof is the observation that the Euler product formula can be generalized. What makes Euler's argument work is the fact that $(1 / m^s) \cdot (1 / n^s) = 1 / (mn)^s$, that is, the fact that the function which maps $n$ to $ 1 / n^s$ is multiplicative. The same argument goes through if we replace the quantity $1 / n^s$ by the function
\[
 \psi(n) = \left\{ \begin{array}{ll}
                      \omega^{\gamma_n} / n^s & \mbox{if $n$ is not divisible by $p$} \\
                      0 & \mbox{otherwise.}
                   \end{array}
            \right.
\]
Thus, generalizing (\ref{euler:product:eqn}), we obtain
\[
 \sum_{p \nmid n} \frac{\omega^{\gamma_n}}{n^s}= \prod_{p \neq q}\left(1-\frac{\omega^{\gamma_q}}{q^{s}}\right)^{-1}.
\]
The sum on the left-hand side ranges over numbers $n$ that are not divisible by $p$, and the product on the right ranges over prime numbers $q$ other than $p$. Euler's calculation then shows that we have
\[
\log\sum_n\frac{\omega^{\gamma_n}}{n^s} =\sum_{q}\frac{\omega^{\gamma_q}}{q^{s}} + O(1), 
\]
in place of (\ref{euler:primes:eqn}). Here the first sum ranges over the same values of $n$, and the second sum ranges over the same values of $p$ as before. \emph{Now} decompose the sum on the right in terms of the remainder that $q$ leaves when divided by $p$, and notice that, by definition, $\gamma_q$ only depends on this remainder. In other words, we have
\begin{multline}
\label{euler:character:eqn}
\log\sum_n\frac{\omega^{\gamma_n}}{n^s} = \Bigg(\sum_{q \equiv 1 \bmod p} \frac{ 1 }{q^{s}}\Bigg) \omega^{\gamma_1} + \Bigg(\sum_{q \equiv 2 \bmod p } \frac{ 1 }{q^{s}}\Bigg) \omega^{\gamma_2} + \ldots + \\
\Bigg(\sum_{q \equiv p - 1 \bmod p } \frac{ 1 }{q^{s}}\Bigg) \omega^{\gamma_{p-1}} + O(1).
\end{multline}

The next step involves the trick we alluded to above. Remember, to show that there are infinitely many primes congruent to $m$ modulo $p$, we want to show that the coefficient of the $m$th term in the preceding equation, $\sum_{q \equiv m \bmod p } \frac{ 1 }{q^{s}}$, approaches infinity as $s$ approaches $1$. If we let $\omega$ be a primitive $(p-1)$st root of 1, then all the roots are given by $\omega^0, \omega^1, \omega^2, \ldots, \omega^{p-2}$. The idea is to plug in all these roots into the preceding equation, and use that to solve for the $m$th coefficient.

Replacing $\omega$ by $\omega^i$ in the last equation yields
\begin{multline*}
\log\sum_n\frac{\omega^{i\gamma_n}}{n^s} = \Bigg(\sum_{q \equiv 1 \bmod p } \frac{ 1 }{q^{s}}\Bigg) \omega^{i \gamma_1} + \Bigg(\sum_{q \equiv 2 \bmod p } \frac{ 1 }{q^{s}}\Bigg) \omega^{i \gamma_2} + \ldots + \\
\Bigg(\sum_{q \equiv p - 1 \bmod p } \frac{ 1 }{q^{s}}\Bigg) \omega^{i \gamma_{p-1}} + O(1).
\end{multline*}
This yields $p-1$ many equations, as $i$ ranges from $0$ to $p - 2$. To solve for the $m$th coefficient, for each $i$, multiply the $i$th equation by $\omega^{-i \gamma_m}$, and add them. 

This is where the magic occurs. If we write $L_i$ for the expression $\sum_{n=1}^{\infty}\frac{\omega^{i\gamma_n}}{n^s}$ that occurs on the left, then the left-hand side of the sum can be written
\[
 \log L_0 + \log L_1 \cdot \omega^{-\gamma_m} + \log L_2 \cdot \omega^{-2 \gamma_m} + \ldots + \log L_{p-2} \cdot \omega^{-(p-2) \gamma_m}.
\]
On the right-hand side, the $m$th term is exactly
\[
(p - 1) \cdot \Bigg(\sum_{q \equiv m \bmod p } \frac{ 1 }{q^{s}}\Bigg),
\]
because $\omega^{i \gamma_m} \cdot \omega^{-i \gamma_m} = 1$ for each $i$, and we are simply summing the same value, $\sum_{q \equiv m \bmod p } 1 / q^{s}$, $p - 1$ times. When $j$ is different from $m$, however, the $j$th term will be
\[
(\omega^{0 (\gamma_j - \gamma_m)} + \omega^{1 (\gamma_j - \gamma_m)} + \ldots + \omega^{(p-2) \cdot (\gamma_j - \gamma_m)}) \cdot \Bigg(\sum_{q \equiv j \bmod p } \frac{ 1 }{q^{s}} \Bigg).
\]
If we write $\eta = \omega^{\gamma_j - \gamma_m}$, then the coefficient in the last expression is
\[
1 + \eta + \eta^2 + \ldots + \eta^{p-2}.
\]
But since $\omega$ is a $(p-1)$st root of 1, so is $\eta$, and since $\gamma_j \neq \gamma_m$, $\eta$ is not equal to $1$. By the observation above, this sum is equal to $0$. In other words, all the other terms magically disappear.

Thus we have shown that
\begin{multline}
\label{main:equation:prime:case:a}
\log L_0 + \omega^{-\gamma_m} \log L_1 + \omega^{-2 \gamma_m} \log L_2 + \ldots + \omega^{-(p-2) \gamma_m} \log L_{p-2} = \\
(p - 1) \cdot \sum_{q \equiv m \bmod p } \frac{ 1 }{q^{s}} + O(1).
\end{multline}
Solving for $\sum_{q \equiv m \bmod p } 1 / q^{s}$ yields
\begin{multline}
\label{main:equation:prime:case}
\sum_{q \equiv m \bmod p } \frac{ 1 }{q^{s}} =
\frac{1}{p - 1} \Bigg(\log L_0 + \omega^{-\gamma_m} \log L_1 + \omega^{-2 \gamma_m} \log L_2 + \ldots + \\
\omega^{-(p-2) \gamma_m} \log L_{p-2} \Bigg) + O(1).
\end{multline}
As a result, we have managed to ``extricate'' the expression $\sum_{q \equiv m \bmod p } 1 / q^{s}$ from (\ref{euler:primes:eqn:b}). The goal is now to show that the this expression approaches infinity as $s$ approaches $1$. We now come to the analytic part of Dirichlet's proof: he showed that as $s$ approaches $1$, $L_0$ approaches infinity, but each of the other $L_i$'s approaches a nonzero limit as $s$ approaches $1$. This implies that the right-hand side approaches infinity as $s$ approaches $1$. Thus the left-hand side approaches infinity as well, which is only possible if there are infinitely many primes congruent to $m$ modulo $p$.

The presentation here follows Dirichlet's short 1837 presentation fairly closely, though Dirichlet is more terse. As Dirichlet pointed out in that note, the argument can be pushed through for an arbitrary modulus $k$. But, as we will see in Section~\ref{dirichlet:section}, the details become unwieldy, and subsequent authors found more convenient ways to express the ideas. In the next section, we explain how the argument above can be described in terms of group characters, and then generalized to the case of an arbitrary modulus.

\subsection{Group characters}
\label{group:character:section}

Let $G$ be a finite abelian group. In contemporary terms, a \emph{character on $G$} is a function $\chi$ from $G$ to the set of nonzero complex numbers with the property that, for every $g_1, g_2 \in G$, $\chi(g_1 g_2) = \chi(g_1) \chi(g_2)$. If $g$ is an element of any finite abelian group, then there is an integer $n > 0$ such that $g^n$ is equal to the identity element of $G$. This implies that $\chi(g)^n = \chi(g^n) = \chi(1) = 1$. This means that for every $g$, $\chi(g)$ is a complex root of 1. The notion of ``character'' introduced in the last section corresponds to the special case where $G$ is the group of nonzero residues modulo $p$, with the operation of multiplication.

The point is that the key properties of the expressions $\omega^{\gamma_n}$ that came into play in the last section hold more generally of the set of characters on any finite abelian group. In particular, for any such group $G$, one can show that there are exactly $|G|$ many distinct characters on $G$, where $|G|$ denotes the number of elements of $G$. In the case where $G$ is the group of nonzero residues modulo $p$, $|G| = p - 1$, so the characters correspond to the $p - 1$ choices of $\omega$ in the previous section. More generally, for any $k \geq 1$, the set of residues $m$ modulo $k$ that have no common factor with $k$ form a group under multiplication. The cardinality of this group is commonly denoted $\ph(k)$, and $\ph$ is known as the Euler phi function. Thus, for every $k$, there are $\ph(k)$ many characters on the group of residues modulo $k$. 

In fact, the set of characters itself has the structure of a group $\widehat G$, where the identity is the character $\chi_1$ that always returns 1, and the product of two characters is given pointwise, $(\chi \cdot \chi')(g) = \chi(g) \chi'(g)$ for every $g$. The following theorem expresses two important properties, known as the ``orthogonality relations'' for group characters.

\begin{theorem}
\label{ortho}
Let $G$ be a finite abelian group.  Then for any character $\chi$ in $\widehat{G}$, we have
\[
\sum_{g \in G}\chi(g) =
\begin{cases}
  |G| & \mbox{if $\chi = \chi_{0}$} \\
  0   & \mbox{if $\chi \neq \chi_{0}$,}
\end{cases}
\]
and for any element $g$ of $G$, we have
\[
\sum_{\chi \in \widehat{G}}\chi(g) =
\begin{cases}
  |G| & \mbox{if $g=1_{G}$} \\
  0   & \mbox{if $g \neq 1_{G}$.}
\end{cases}
\]
\end{theorem}

The remarkable fact is that it is no harder to prove these facts in the general case than in the specific case where $G$ is a group of residues modulo $p$. For example, the second equation clearly holds when $g$ is the identity of $G$, since, in this case, each term of the sum is equal to $1$. Otherwise, choose a character $\psi$ such that $\psi(g) \neq 1$ and note
\[
 \psi(g) \sum_{\chi \in \widehat{G}} \chi(g) = \sum_{\chi \in \widehat{G}} \psi(g) \chi(g) = \sum_{\chi \in \widehat{G}} \chi(g),
\]
since multiplying each character $\chi$ in $\widehat{G}$ by $\psi$ simply permutes the elements of $\widehat{G}$. Subtracting the right side of the equation from the left, we see that $(\psi(g) - 1) \cdot \sum_{\chi \in \widehat{G}} \psi(g) = 0$, and since $\psi(g)$ is not equal to $1$, we have that $\sum_{\chi \in \widehat{G}} \psi(g) = 0$. The first equation can be established in a similar way.

The second orthogonality relation gives rise to the ``cancellation trick'' used in the last section, where we multiplied each identity by $\omega^{- i \gamma_m}$ and added them, to isolate a particular coefficient. The general phenomenon can be expressed as follows:
\begin{corollary}
\label{orthocorrol}
For any $g, h \in G$ we have the following:
\[
\sum_{\chi \in \widehat{G}}\chi(g)\overline{\chi(h)}= \begin{cases} |G| & \ \mbox{\rm if} \ g=h \\ 0 & \ \mbox{\rm if} \ g \neq h.\end{cases}
\]
\end{corollary}
Here $\bar z$ denotes the complex conjugate of $z$, which is in fact equal to $1 / z$ when $z$ is a root of unity. The corollary follows from the fact that we have
\[
\sum_{\chi \in \widehat{G}}\chi(g)\overline{\chi(h)} =\sum_{\chi \in \widehat{G}}\chi(g)\chi(h)^{-1}=\sum_{\chi \in \widehat{G}}\chi(gh^{-1})= \begin{cases} |G| & \ \mbox{if} \ g=h \\ 0 & \ \mbox{if} \ g \neq h. \end{cases}
\]

Notice that the abstract algebraic formulation simplifies matters by eliminating clutter. For example, the presentation in the last section depended on choices of a primitive element $g$ modulo $p$, and a primitive $(p-1)$st root of unity $\omega$. Although these played a role in the computations, any choice of $g$ and $\omega$ works just as well. The abstract version ``factors these out'' of the presentation. Recall also that the calculation in the last section required facts such as $\gamma_{m n} = \gamma_m + \gamma_n$. Once again, the abstract version factors this out of the computation; the requisite property of $\gamma$ subsumed by the more general fact that $\widehat{G}$ is a group, and only the latter fact enters into the proof.

\subsection{A modern formulation of Dirichlet's proof}
\label{modern:formulation:section}

With the notion of a group character in mind, we can now describe Dirichlet's original proof of Theorem~\ref{dirichlet:theorem} in modern terms. Let $k$ be an integer greater than or equal to $1$. It is a fundamental theorem of number theory that an integer $n$ is relatively prime to $k$ if and only if $n$ has a multiplicative inverse modulo $k$; in other words, if and only if there is some $n'$ such that $n n' \equiv 1 \bmod{k}$. This implies that the residues of integers modulo $k$ that are relatively prime to $k$ form a group, denoted $(\ZZ/k\ZZ)^*$, with multiplication modulo $k$. As noted above, the cardinality of $(\ZZ/k\ZZ)^*$, that is, the number of residues relatively prime to $k$, is denoted $\ph(k)$.

A character $\chi$ on the group of residues modulo $k$ can be viewed as a function defined on all integers by
\[
 X(n) =
  \begin{cases}
    \chi(n \bmod k) &\text{if $n$ is relatively prime to $k$} \\
    0 &\text{otherwise.}
  \end{cases}
\]
Such a function is called a \emph{Dirichlet character modulo $k$}. Dirichlet characters are \emph{completely multiplicative}, which is to say, $X(1) = 1$ and $X(mn) = X(m)X(n)$ for every $m$ and $n$ in $\ZZ$. Mathematicians typically use the symbol $\chi$ to range over Dirichlet characters, blurring the distinction between such functions and their group-character counterparts. This is harmless, since there is a one-to-one correspondence between the two, and so we will adopt this practice as well.

Recall that in the case where $k$ is a prime number $p$, Dirichlet considered certain expressions $L_i(s)$, analogues of Euler's zeta function, where $i$ is an integer between $0$ and $p - 2$. Each such $i$ corresponds to a choice of a character $\chi$ modulo $p$. In the modern formulation, then, we define
\[
  L(s, \chi)= \sum_{n=1}^{\infty}\frac{\chi(n)}{n^{s}},
\]
where $\chi$ is such a character. The function $L(s, \chi)$ is called the \emph{Dirichlet $L$-function}, or \emph{$L$-series}.

The calculation in Section~\ref{dirichlet:approach:section} can be generalized to show:
\[
\log L(s, \chi) = \sum_{q \nmid k}\frac{\chi(q)}{q^{s}} \ + \ O(1).
\]
Now comes the crucial use of Corollary~\ref{orthocorrol} to pick out the primes in the relevant residue class.  We multiply each side of the above equation by $\overline{\chi(m)}$ and then take the sum of these over all the Dirichlet characters modulo $k$. (Recall that we can identify each Dirichlet character with the corresponding group character, that is, the corresponding element of $\widehat{(\mathbb{Z}/k \mathbb{Z})^*}$.) Thus we have:
\begin{equation*}
% \label{sumoverchar1}
\sum_{\chi \in \widehat{(\mathbb{Z}/k \mathbb{Z})^*}} \overline{\chi(m)}\log L(s, \chi) =
\sum_{\chi \in \widehat{(\mathbb{Z}/k \mathbb{Z})^*}} \overline{\chi(m)}\sum_{q \nmid k}\frac{\chi(q)}{q^{s}}\  + \  O(1).
\end{equation*}
To simplify this expression, we exchange the summations on the right-hand side, and appeal to Corollary~\ref{orthocorrol}. Since the cardinality of the group $(\mathbb{Z} / k \mathbb{Z})^*$ is $\ph(k)$, we obtain
\begin{equation}
\label{sumoverchar2}
\sum_{\chi \in \widehat{(\mathbb{Z} / k \mathbb{Z})^*}}\overline{\chi(m)}\log L(s, \chi) =
\ph(k) \sum_{q \equiv m \pmod{k}}\frac{1}{q^{s}} \ + \ O(1).
\end{equation}
This is analogous to the equation (\ref{euler:primes:eqn}) in Euler's proof, and equation (\ref{main:equation:prime:case:a}) in Section~\ref{dirichlet:approach:section}. Our goal is once again to show that the left-hand side tends to infinity as $s$ approaches 1 from above. This implies that the right-hand side tends to infinity, which, in turn, implies that there infinitely many primes $q$ that are congruent to $m$ modulo $k$.

To show that $\sum_{\chi \in \widehat{(\mathbb{Z} / k \mathbb{Z})^*}}\overline{\chi(m)}\log L(s, \chi)$ tends to infinity as $s$ approaches $1$, we divide the characters into three classes, as follows:
\begin{enumerate}
\item The first class contains only the principal character $\chi_0$, which takes the value of 1 for all arguments that are relatively prime to $k$, and 0 otherwise.
\item  The second class consists of all those characters which take only real values (i.e.\ 0 or $\pm 1$), other than the principal character.
\item  The third class consists of those characters which take at least one complex value.
\end{enumerate}
It is not difficult to show that $L(s, \chi_{0})$ has a simple pole at $s=1$, which implies that the term $\overline{\chi_0}(m) \log L(s,\chi_0)$ approaches infinity as $s$ approaches $1$. The real work involves showing that for all the other characters $\chi$, $L(s,\chi)$ has a finite nonzero limit. This implies that the other terms in the sum approach a finite limit, and so the entire sum approaches infinity.

For characters in the third class, that is, the characters that take on at least one complex value, the result is not difficult. For characters in the second class, the result is much harder, and Dirichlet used deep techniques from the theory of quadratic forms to obtain it. In the years that followed, other mathematicians found alternative, and simpler, ways of handling this case. But even in modern presentations, this case remains the most substantial and technically involved part of the proof.

\section{Functions as objects}
\label{functions:section}

In Section~\ref{dirichlet:section} below, we will discuss, in greater detail, the implicit treatment of characters in Dirichlet's original proof, and in Section~\ref{transition:section}, we will summarize the gradual historical transition to the modern formulation. The general theme will be that, over time, characters came to be treated as objects in their own right. Before surveying the history, however, it will be helpful for us to provide some general background information on the nineteenth century concept of ``function,'' and begin to spell out what it means to treat functions like characters as ``objects.''

In ``Concept,'' we discussed a number of nineteenth-century methodological changes that are clustered around the function concept. These include what we termed the ``unification'' or ``generalization'' of the function concept, whereby particular instances (including real- and complex-valued functions, number-theoretic functions, sequences, permutations, transformations, automorphisms, and so on) gradually came to be subsumed under a general notion; the ``liberalization'' of the function concept, whereby mathematicians adopted novel means of defining particular functions, such as Dirichlet's 1827 example of a real-valued function that takes one value on the rationals, and another value on the irrationals; the ``extensionalization'' of the function concept, whereby functions gradually came to be viewed less as syntactic or algebraic expressions, and more as the abstract entities denoted by such expressions; and the ``reification'' of the function concept, whereby functions were gradually treated as \emph{bona fide} mathematical objects.

The notion of ``reification'' is vague. The claim that over the course of the century characters gradually become treated as new sorts of objects supports our contention that the transformation has ontological overtones, but it raises serious questions as to what, exactly, it means to treat certain entities as objects.

To start with, consider the fact that in our presentation of Dirichlet's theorem we identified the concept of a ``character,'' reasoned about the entities falling under this concept, and ascribed various properties to them. This seems to be a bare-minimum requirement to support the claim that a mathematical text sanctions certain entities as objects, namely, that it recognizes them as being entities of a certain \emph{sort}, capable of bearing predicates and being the target of certain operations. It does not matter whether we take this sort as fundamental (for example, as we take the notion of ``integer'' in most contexts) or as derived from a broader sort (for example, when we view characters as functions of a certain kind). What is important is that the entities belong to a grammatically recognized category, and this category helps determine the predicates and operations that can be meaningfully ascribed to it. For example, one can talk about one integer being larger than another, but not one character as being larger than another. In sum, our first criterion of objecthood is whether the entities in question have a recognizable role in the grammar of the language.

The fact that we took characters to be ``represented'' by certain symbolic expression provides another clue, insofar as we generally speak of a representation \emph{of} something or other. For example, we think of expressions like ``$6$'' and ``$2 \times 3$'' as representing an integer. As Michael Detlefsen has pointed out to us, one common view is that an ``object'' is what remains invariant under all its representations; in other words, what is left over when one has ``squeezed out'' all the features that are contingent on particular representations. When it comes to the notion of a function, what is the underlying invariant? There may be lots of ways of describing a particular function, but what makes them representations of the \emph{same} function is surely that they take the same values on any given input. Thus treating function expressions extensionally is a sign that one is reasoning about functions as objects, rather than reasoning about the expressions themselves.\footnote{Recall Quine's dictum that ``there is no entity without identity,'' for example in \cite{quine:69}.}

A third hallmark of object-hood that is present in our list is evidenced by the fact that we can \emph{sum} over characters, just as we can sum over natural numbers. Notice that in an expression $\sum_{\chi} \ldots \chi \ldots$, the variable $\chi$ is a bound variable that ranges over the entities in question. Similar considerations hold for the universal and existential quantifiers. Viewing the natural numbers as quintessential mathematical objects, a sign that an entity has attained the status of object-hood is that it is possible to quantify over them in theorems and definitions, just as one quantifies over the natural numbers.\footnote{This echoes another Quine dictum, ``to be is to be the value of a bound variable'' \cite[p.~15]{quine:48}.} The consideration admits of degrees: whereas the bare-minimum requirement discussed above may allow us to state theorems about, and define operations on, ``arbitrary'' entities of the sort, a more
full-blown notion of object-hood will give us more latitude in the kinds of quantification and binding that are allowed.

A fourth criterion for object-hood is evidenced by the fact that characters are allowed to appear as \emph{arguments} to the $L$-functions, for example, in the expression $L(s,\chi)$. To avoid making this consideration depend on the modern notion of a function, let us note that what is essential here is that an expression denoting a recognized mathematical object (in this case, a complex number) is allowed to \emph{depend} on a character, much the way that a real number $(s)_i$ in a sequence depends on the value of the index $i$, or a value $\ph(n)$ of the phi function depends on $n$. What makes this more potent than the mere ability to define operations on characters is that the dependent expressions are treated as objects in their own right. $L(s,\chi)$ is not just an operation on $s$ and $\chi$: fixing $\chi$, the function $s \mapsto L(s,\chi)$ is an object that one can integrate and differentiate, and fixing $s$, we can sum over the values obtained by varying $\chi$. 

It is also notable that the characters can be components in the construction of other mathematical objects and structures. For example, one can form sets and sequences of characters, in much the same way that one forms sets and sequences of numbers, and one can define a group whose elements are characters, in much that same way that one can form a group whose elements are residues modulo some number $m$.

To summarize, here are some of the various senses in which one might say that characters are treated ``as objects'' in our presentation of Dirichlet's proof:
\begin{enumerate}
 \item Characters fall under a recognized grammatic category, which allows us to state things about them and define operations and predicates on them.
 \item There is a clear understanding of what it means for two expressions to represent the \emph{same} character, and conventions ensure that the expressions occurring in a proof respect this ``sameness.''
 \item One can quantify and sum over characters; in short, they can fall under the range of a bound variable.
 \item One can define functions which take characters as arguments.
 \item One can construct new mathematical entities, like sets and sequences, whose elements are characters. In particular, characters can be elements of an algebraic structure like a group. 
\end{enumerate}
We recognize that determining the ``ontological commitments'' of a practice may not be as clear-cut as Quine's writings suggest. Our goal here is not to explicate what it means to say that a certain manner of discourse is committed to treating some entity as an object. In particular, we do not claim to have given a precise sense to the question as to whether a particular mathematical proof is committed to functions as objects. We do claim, however, to have identified various important senses in which contemporary proofs of Dirichlet's theorem treat functions as ordinary mathematical objects, whereas Dirichlet's original proof did not.

It may be helpful to compare the way we treat functions today to the way we treat natural numbers today. For example, the expressions ``$2 + 2$'' and ``$4$'' both denote integers, but we think of the number as the object denoted, rather than the expression. Thus we can send numbers as arguments to functions, and when we write $f(2 + 2)$ and $f(4)$, it is understood that the function $f$ cannot distinguish the mode of presentation. We can form sets of numbers, like the set of even number or the set of prime numbers, and we can consider algebraic structures on these sets; for example, the ring of integers, or the field of integers modulo $7$. We can quantify over numbers in definitions, such as when we say $n$ divides $m$ if there is some $k$ such that $n k = m$, and in theorems, such as when we assert that every integer greater than one has a prime divisor. If $S$ is a finite set of integers and $f$ is a function from the integers to the integers or the reals, we can readily form the sum $\sum_{x \in S} f(x)$.

In contemporary mathematics, nothing goes awry if you replace integers with functions in the examples in the last paragraph. In other words, one can define functionals $F(f)$ that depend only on the extension of $f$, and not its manner of presentation. We can consider sets of functions, rings of functions, and spaces of functions. We quantify over functions in definitions and theorems, and, if $S$ is a finite set of functions, we think nothing of considering a sum $\sum_{f \in S} F(f)$. In the proof of Dirichlet's theorem, these ``higher order'' operations are manifest when we consider the group of characters $\chi$, define the Dirichlet $L$ series $L(s, \chi)$, and form the sum $\sum_\chi \overline{\chi(m)}\log L(s, \chi)$.

In ``Concept,'' we argued in detail that these very features of the modern treatment of functions were alien to early nineteenth century mathematics, and that the history of presentations of Dirichlet's theorem shows a very gradual evolution, in fits and starts, towards the contemporary manner of thought. We will highlight some of the key features of the historical development in Sections~\ref{dirichlet:section} and \ref{transition:section}, and, in Section~\ref{analysis:section}, explore what the history tells us about the nature of mathematics.

\section{Dirichlet's treatment of characters}
\label{dirichlet:section}

Contemporary mathematicians are often surprised to hear that there is no explicit notion of ``character'' in Dirichlet's 1837 proof. After all, the expressions $X(n)$ defined in Section~\ref{modern:formulation:section} are known as ``Dirichlet characters'' precisely because of their implicit use in that proof. But Dirichlet did not introduce notation for the characters or refer to them as such. When we speak of the ``characters'' in his proof, we are projecting a modern interpretation onto the symbolic expressions that appear there.

Remember how it works in the case where the common difference is a prime, $p$. Let $g$ be a primitive element modulo $p$, and for every $n$ coprime to $p$, let $\gamma_n$ denote the index of $n$ with respect to $g$, so that $g^{\gamma_n} \equiv n \bmod p$. Then each character $\chi$ corresponds to a $(p-1)$st root of unity $\omega$, with defining equation $\chi(n) = \omega^{\gamma_n}$. In that case, Dirichlet wrote $\omega^{\gamma_n}$ where we would write $\chi(n)$.

We obtain all the characters by picking a primitive $(p-1)$st root of unity, $\Omega$, so that all the $(p-1)$st roots of unity are given by the sequence $\Omega^0, \ldots, \Omega^{p-2}$. This provides a convenient numbering scheme for the characters and $L$-series: Dirichlet used $L_m$ to denote the $L$-series based on the character $\chi$ that corresponds to $\Omega_m$, where we would write instead $L(s, \chi)$. And where we would form a summation over the set of all characters, Dirichlet instead took a summation over the values $0, \ldots, p-2$. For example, after demonstrating the Euler product formula,
\[
\prod \frac{1}{1-\omega^{\gamma}\frac{1}{q^{s}}}=\sum\omega^{\gamma}\frac{1}{n^{s}}=L,
\]
Dirichlet wrote: 
\begin{quote}
The equation just found represents $p-1$ different equations that result if we put for $\omega$ its $p-1$ values.  It is known that these $p-1$ different values can be written using powers of the same $\Omega$ when it is chosen correctly, to wit:
\[ 
\Omega^{0}, \ \Omega^{1}, \ \Omega^{2}, \  \ldots,\ \Omega^{p-2}
\]
According to this notation, we will write the different values $L$ of the series or product as:
\[
L_0,\ L_1,\ L_2,\ \ldots,\ L_{p-2}
\]
\end{quote}

In the case where the modulus $k$ is not prime, the procedure is more complicated. It is a fundamental theorem of group theory that every finite abelian group can be represented as a product of cyclic groups, but that theorem was first proved by Kronecker in 1870 \cite{kronecker:70}. Dirichlet instead used the particular instance of this fact for the group $(\ZZ / k\ZZ)^*$ of residues modulo $k$ that are relatively prime with $k$ (these are sometimes called the ``units'' modulo $k$). The structure of that group was known to Gauss. First, write $k$ as a product of primes,
\[
 k = 2^\lambda p_1^{\pi_1} p_2^{\pi_2} \cdots p_j^{\pi_j}
\]
where each $p_i$ is an odd prime and $\pi_i$ is greater than or equal to 1. Then the group of units modulo $k$ is isomorphic to the product of the groups of units modulo each term in the factor. If $p$ is an odd prime and $\pi$ is an integer greater than or equal to $1$, then one can more generally find a primitive element $c$ modulo $p^\pi$. This means that the residue class of $c$ generates the cyclic group $(\ZZ / p^\pi \ZZ)^*$, or, equivalently, for every $n$ relatively prime to $p$ there is a $\gamma_n$ such that $c^{\gamma_n} \equiv n \bmod p^\pi$. Thus we can choose primitive elements $c_1, \ldots, c_j$ corresponding to $p_1^{\pi_1}, p_2^{\pi_2}, \ldots, p_j^{\pi_j}$. If $\lambda \geq 3$, however, there is no primitive element modulo $2^\lambda$. Rather, $(\ZZ/ 2^\lambda)^*$ is a product of two cyclic groups, and for every $n$ relatively prime to $2^\lambda$ there are an $\alpha_n$ and $\beta_n$ such that $(-1)^{\alpha_n} 5^{\beta_n} \equiv n \bmod 2^\lambda$. Thus for any $n$ relatively prime to $k$, we can write
\[
 n \equiv (-1)^{\alpha_n} 5^{\beta_n} c_1^{\gamma_{1,m}} c_2^{\gamma_{2,m}} \ldots c_j^{\gamma_{j,m}} \bmod k
\]
where each $\gamma_{i,n}$ is the index $n$ relative to $p_i^{\pi_i}$. As above, if we choose appropriate roots of unity $\theta, \ph, \omega_1, \omega_2, \ldots, \omega_j$, we obtain a character
\begin{equation*}
%\label{characterconstruction}
 \chi(n) = \theta^{\alpha_n} \ph^{\beta_n} \omega_1^{\gamma_{1,n}} \omega_2^{\gamma_{2,n}} \cdots
 \omega_j^{\gamma_{j,n}}.
\end{equation*}
And, once again, every character is obtained in this way. We should note that Dirichlet used the notation $p, p', \ldots$ rather than $p_1, \ldots, p_j$ to denote the sequence of odd primes. Moreover, he used the notation $\alpha, \beta, \gamma, \gamma', \ldots$ to denote the indices, suppressing the dependence on $n$. Thus, Dirichlet wrote $\theta^\alpha \ph^\beta \omega^\gamma \omega'^{\gamma'} \ldots$ for the expression we have denoted $\chi(n)$ above, leaving it up to us to keep in mind that $\alpha, \beta, \ldots$ depend on $n$.

To summarize, in the simple case of a prime modulus $p$, Dirichlet fixed a primitive element modulo $c$, and represented each character $\chi$ in terms of a $(p-1)$st root of unity, $\omega$. In that case, the value $\chi(n)$ is given by $\omega^{\gamma_n}$. In the more general case of a composite modulus $k$, Dirichlet fixed primitive elements modulo the terms of the prime factorization of $k$, and represented each character $\chi$ in terms of a sequence $\theta, \ph, \pi, \pi'$ of roots of unity. In that case, the value $\chi(n)$ was written $\theta^\alpha \ph^\beta \omega^\gamma \omega'^{\gamma'} \ldots$, suppressing the information that the exponents $\alpha, \beta, \gamma, \gamma', \ldots$ depend on $n$. For example, he described the Euler product formula as follows:

\begin{quote}
\begin{align}
 \prod \frac{1}{1-\theta^{\alpha}\ph^{\beta}\omega^{\gamma}\omega^{'\gamma^{'}}\ldots\frac{1}{q^{s}}}=\sum \theta^{\alpha}\ph^{\beta}\omega^{\gamma}\omega'^{\gamma'}\ldots\frac{1}{n^{s}} = L, \label{dirichletEulergeneral}
\end{align}
where the multiplication sign ranges over all primes, with the exclusion of $2, p, p', \ldots$, and the summation ranges over all the positive integers that are not divisible by any of the primes $2, p, p', \ldots$. The system of indices $\alpha, \beta, \gamma, \gamma', \ldots$ on the left side corresponds to the number $q$, and on the right side to the number $n$. The general equation (\ref{dirichletEulergeneral}), in which the different roots $\theta, \phi, \omega, \omega', \ldots$ can be combined with one another arbitrarily, clearly contains $K$-many particular equations. \cite[p.~17; equation number changed]{Dirichlet37}
\end{quote}
Note, again, Dirichlet's characterization of the general equation as ``containing'' the particular instances.
Here, $K$ is what we have called $\ph(k)$, the cardinality of the group $(\ZZ / k\ZZ)^*$. 

Dirichlet went on to observe that we can choose primitive roots of unity $\Theta, \Phi, \Omega, \Omega', \ldots$ so that all choices of $\theta, \ph, \omega, \omega', \ldots$ can be expressed as powers of these,
\[
\theta = \Theta^\mathfrak{a}, \ph = \Phi^\mathfrak{b}, \omega = \Omega^\mathfrak{c}, \omega' = \Omega^{\mathfrak{c}'}, \ldots,
\]
just as in the simpler case. He wrote that we can thus refer to the $L$-series in a ``convenient'' (\emph{bequem}) way, as $L_{\mathfrak{a}, \mathfrak{b}, \mathfrak{c}, \mathfrak{c}', \ldots}$, where $\mathfrak{a}, \mathfrak{b}, \mathfrak{c}, \mathfrak{c}', \ldots$ are the exponents of the chosen primitive roots. Notice that the representations just described depend on fixed, but arbitrary, choices of the primitive roots of unity, as well as fixed but arbitrary generators of the cyclic groups. Modulo those choices, we have parameters $\mathfrak{a}, \mathfrak{b}, \mathfrak{c}, \mathfrak{c}', \ldots$ that vary to give us all the characters; and for each choice of $\mathfrak{a}, \mathfrak{b}, \mathfrak{c}, \mathfrak{c}', \ldots$ we have an explicit expression that tells us the value of the character at $n$. 

For Dirichlet, summing over characters therefore amounted to summing over all possible choices of this representing data. In the special case of where the common difference is a prime, $p$, Dirichlet ran through calculations similar to those described in Section~\ref{modern:formulation:section} to obtain the following identity:
\begin{multline*}
% \label{Dirichletsum1}
\sum\frac{1}{q^{1 + \rho}} +\frac{1}{2}\sum\frac{1}{q^{2 +2 \rho}} + \frac{1}{3}\sum\frac{1}{q^{3 + 3 \rho}} + \ldots \\
= \frac{1}{p-1} (\log L_{0} + \Omega^{-\gamma_{m}} \log L_{1} + \Omega^{-2\gamma_{m}}\log L_{2} + \ldots + \Omega^{-(p-1)\gamma_{m}} \log L_{p-2}).
\end{multline*}
This is exactly equation (\ref{main:equation:prime:case}) above, with $\Omega$ in place of our $\omega$, and $1 + \rho$ in place of $s$, and the ``$O(1)$'' expression left explicit. In the more general case, he arrived at the analogous result:
\begin{multline*}
%  \label{Dirichletsum2}
\sum\frac{1}{q^{1 + \rho}} + \frac{1}{2}\sum\frac{1}{q^{2 +2 \rho}} + \frac{1}{3}\sum\frac{1}{q^{3 + 3 \rho}} + \ldots \notag \\
= \frac{1}{K}\sum \Theta^{-\alpha_{m}\mathfrak{a}}\ \Phi^{-\beta_{m}\mathfrak{b}}\Omega^{-\gamma_{m}\mathfrak{c}}\Omega^{-\gamma'_{m}\mathfrak{c}'} \cdots \log L_{\mathfrak{a},\mathfrak{b},\mathfrak{c},\mathfrak{c}', \ldots.}
\end{multline*}
Here the summation on the right-hand side of the equation is over the possible values of $\mathfrak{a}, \mathfrak{b}, \mathfrak{c}, \mathfrak{c'}, \ldots$. This corresponds to equation (\ref{sumoverchar2}) in Section~\ref{modern:formulation:section}.

Finally, recall from the sketch in Section~\ref{modern:formulation:section} that Dirichlet divided the $L$ functions into three classes, depending on whether the corresponding character was trivial (identically equal to 1), real-valued, or complex-valued. But in Dirichlet's presentation, the categorization was made in terms of the \emph{roots used to describe the character}. Thus the three classes of $L$ functions were characterized as follows:
\begin{enumerate}
 \item the one in which all the roots contained in the expression are $1$
 \item those, among the ones that remain, in which all the roots are real ($\pm 1$)
 \item those in which at least one of the roots is not real
\end{enumerate}
Dirichlet showed that the first approaches infinity as $\rho$ approach $0$, while the others approach finite limits, which establishes the desired conclusion.

Let us summarize the features of Dirichlet's presentation we wish to highlight. First, he did not name or identify the characters, and simply used the corresponding algebraic expressions. The corresponding $L$ functions were then characterized by the data that appeared in the expression, rather than in terms of a functional dependence on the character. In other words, Dirichlet wrote $L_m$ or $L_{a,b,c,c', \ldots}$ where we would write $L(s, \chi)$. As a result, where we would sum an expression over all values of the characters $\sum_{\chi} \ldots$, he summed over the representing data $\sum_m \ldots$ or $\sum_{a, b, c, c', \ldots} \ldots$. Finally, in preparation for the analytic part of the proof, he sorted the $L$ functions in terms of this data, rather than in terms of the values of the corresponding characters. 

In the next section, we will see that, over time, all of these features were gradually eliminated from later expositions.

\section{The transition to the modern treatment of characters}
\label{transition:section}

In ``Concept,'' we studied the treatment of characters in subsequent work by Dirichlet (1840, 1841), Dedekind (1863, 1879), Kronecker (1870's), Weber (1883), Hadamard (1896), de la Vall\'ee Poussin (1897), and Landau (1909, 1927). We will not review all the details here, but, rather, summarize the salient features of the history. 

% In 1840 and 1841, Dirichlet sketched extensions of the result to the Gaussian integers and to the set of numbers represented by certain quadratic forms, respectively. In 1863, in an appendix to his write-up of Dirichlet's lectures on number theory, Dedekind presented a proof of Dirichlet's theorem on primes in an arithmetic progression as well. In 1879, Dedekind introduced the notion of a character on an abelian group other than the integers modulo some number $k$, and in 1883, Weber provided a sustained study of the properties of characters on an arbitrary finite abelian group. Hadamard, in 1896, and de la Vall\'ee Poussin, in 1897, presented the central methods of Dirichlet's proof, towards generalizing them. In 1901, Kronecker's student Hensel gave a presentation of a constructive proof that Kronecker had given in the 1870's and 1880's, in Hensel's write-up of Kronecker's lectures on number theory. Finally, Landau gave a presentation of Dirichlet's proof in 1909, and again, in 1927, in an elementary textbook on number theory. 

\subsection{Reification}

We have seen that in Dirichlet's original proof, characters are present only in the form of the algebraic expressions $\omega^{\gamma_n}$ in the simple case, and in the form $\theta^{\alpha_n} \ph^{\beta_n} \omega^{\gamma_n} \omega'^{\gamma'_n} \ldots$ in the case of an arbitrary modulus. In 1841, however, Dirichlet considered expressions
\[
\Omega_{n}=\varphi^{\alpha_{n}}\varphi^{'\alpha'_{n}}\times\ldots\times\psi^{\beta_{n}}\chi^{\gamma_{n}}\psi^{'\beta'_{n}}\chi^{'\gamma'_{n}}\times\ldots\times\theta^{\delta_{n}}\eta^{\varepsilon_{n}}
\]
analogous to the characters in his 1837 proof. In this case, however, he introduced the explicit notation $\Omega_n$, and isolated four key properties of these values:
\begin{enumerate}
\item $\Omega_{nn'} = \Omega_{n}\Omega_{n'}$ for every $n$ and $n'$.
\item $\Omega_{n'}=\Omega_{n}$ whenever $n'\equiv n \pmod{k}$.
\item $\sum\Omega_{l} = 0$ or $\sum\Omega_{l} = \frac{1}{4}\psi(k)$ depending on whether there is at least one root among the roots in $\Omega_{l}$ that is different to 1, or whether they are all equal.
\item $S\Omega_{n} = \frac{1}{4}\psi(k)$ or $S\Omega_{n}=0$ depending on whether $n\equiv 1 \pmod{k}$ or $n\not\equiv 1 \pmod{k}$, where the sign ``$S$'' indicates a sum over all combinations of the roots that can occur in $\Omega$.
\end{enumerate}
In modern terms, the first clause asserts that the function $n \mapsto \Omega_n$ is a multiplicative function from the integers to the complex numbers, and the second asserts that the value $\Omega_n$ only depends on the value of $n$ modulo $p$. If you add the constraints that $\Omega_n$ is nonzero when $n$ is relatively prime to $k$ and zero otherwise, this is exactly the algebraic definition of character we presented in Section~\ref{group:character:section}. The third and fourth properties correspond to the two orthogonality relations we presented in Section~\ref{group:character:section}. The article provided only a short sketch of a generalization of his 1837 proof, but it is notable that there  Dirichlet went out of his way to flag these expressions as playing a key role, and to abstract away the general properties that are common to both proofs. 

In 1863, Dedekind gave an exposition of Dirichlet's proofs in one of the appendices, or ``supplements,'' to the first edition of his presentation of Dirichlet's lectures on number theory \cite{dirichlet:63b}. When presenting the generalization of the Euler product formula, he went out of his way to point out that the function
\[
 \psi(n) = \frac{\theta^{\alpha}\eta^{\beta}\omega^{\gamma}\omega'^{\gamma'}\ldots}{n^s}
\]
is multiplicative, and that this is what makes the generalization hold. In a later 1871 edition of the work, he added a footnote, in which he singled out the numerator of this expression, and introduced the notation $\chi(n)$:
\begin{quote}
The numerator [of $\psi(n)$] $\chi(n)=\theta^{\alpha}\eta^{\beta}\omega^{\gamma}\omega'^{\gamma'}\ldots$ has the characteristic property $\chi(n)\chi(n')=\chi(nn')\ldots$\cite[\S 133, footnote]{dirichlet:63b}
\end{quote}
It is notable that he went out of his way to add this footnote, calling attention to the importance of these expressions.\footnote{In ``Concept,'' we mistakenly asserted that Dedekind did not alter the text of this supplement in later editions. He made very few such changes, however, making this particular addition especially interesting.}

In 1879, in the third edition of the lectures, Dedekind introduced the notion of a character in an entirely different context: his theory of ideals in an algebraic number field. Rather than considering characters on the multiplicative group of residues modulo an integer, he considered characters defined on another finite abelian group, namely, on the class group in an algebraic number field: 
\begin{quote}
 \ldots the function $\chi(\mathfrak a)$ also possesses the property that it takes the same value on all ideals $\mathfrak a$ belonging to the same class $A$; this value is therefore appropriately denoted by $\chi(A)$ and is clearly always an $h$th root of unity. Such functions $\chi$, which in an extended sense can be termed \emph{characters}, always exist; and indeed it follows easily from the theorems mentioned at the conclusion of \S 149 that the class number $h$ is also the number of all distinct characters $\chi_1, \chi_2, \ldots, \chi_h$ and that every class $A$ is completely characterized, i.e.~is distinguished from all other classes, by the $h$ values $\chi_1(A), \chi_2(A), \ldots, \chi_h(A)$.\footnote{The quotation appears in \S 178 in the 1879 edition of the \emph{Vorlesungen} \cite{dirichlet:63b}, and in \S 184 of the 1894 edition, which is reproduced in Dedekind's \emph{Werke} \cite{dedekind:68}. The translation above is by Hawkins \cite[p.~149]{hawkins:71}.}
\end{quote}
As we emphasize in ``Concept,'' this was not only the first use of the term ``character'' in its modern sense, but also, as far as we know, the earliest instance of the use of the term ``function'' for something defined on a domain other than the integers, real numbers, or complex numbers. (A similarly broad use of the term occurs in Frege's \emph{Begriffsschrift}, which was published in the same year.) We will discuss Frege's notion of function in detail in Sections~\ref{frege:section} and \ref{frege:section:b}.) Within three years, in an 1882 publication, Weber gave the general definition of a character of an abelian group and provided a thorough analysis of their properties.

Thus, over time, the symbolic expressions appearing in Dirichlet's proof were named and flagged as entities worthy of attention. Their properties were stated abstractly, and developed in a manner that were independent of the original formulation. This, in turn, made it possible to apply the notion in other settings. As we have argued in Section~\ref{functions:section}, this provides at least a minimal sense in which characters can be viewed as objects, namely, as entities which can bear properties and be a target of assertions.

\subsection{Functional dependence and summation}

In Section~\ref{functions:section}, we also flagged it as notable that, in the modern view, functions can depend on characters, and we can form the sum of an expression with a variable ranging over the characters. Let us consider the way these features of the treatment of characters play out in the various presentations of Dirichlet's theorem.

We have noted that one benefit of identifying the characters as such is that it facilitates extracting the central properties that play a role in the proof, such as the identity
\[
\sum_{\chi \in \widehat{G}}\chi(g) =
\begin{cases}
  |G| & \mbox{if $g=1_{G}$} \\
  0   & \mbox{if $g \neq 1_{G}$}
\end{cases}
\]
in Theorem~\ref{ortho}, and the consequence expressed by Corollary~\ref{orthocorrol} that for every $g$ and $h$ in an abelian group $G$, 
\[
\sum_{\chi \ \in \ \widehat{G}}\chi(g)\overline{\chi(h)}= \begin{cases} |G| & \ \mbox{if} \ g=h \\ 0 & \ \mbox{if} \ g \neq h.\end{cases}
\]
In the case where $G$ is the group of nonzero residues modulo $p$, Dirichlet expressed the latter by saying that we have
\[
1+\Omega^{h\gamma - \gamma_{m}} + \Omega^{2(h\gamma - \gamma_{m})} + \ldots + \Omega^{(p-2)(h\gamma - \gamma_{m})}=0 
\]
except when $h\gamma - \gamma_{m} \equiv 0 \bmod{p-1}$, in which case the sum is equal to $p - 1$. In the case of an arbitrary modulus, Dirichlet did not even extract the conclusion explicitly. Rather, it is implicitly contained in an argument in which he considered the sum $\frac{1}{h}\sum W\frac{1}{q^{h + h\rho}}$,
\begin{quote}
\ldots where the symbol $\sum$ ranges over all primes $q$ and $W$ denotes the product of the sums taken over $\mathfrak{a}, \mathfrak{b}, \mathfrak{c}, \mathfrak{c'}, \ldots$ or respectively over 
\[                                                                                                       
\sum\Theta^{(h\alpha - \alpha_{m})\mathfrak{a}}, \sum\Phi^{(h\beta - \beta_{m})\mathfrak{b}}, \sum\Omega^{(h\gamma - \gamma_{m})\mathfrak{c}}, \sum\Omega'^{(h\gamma' - \gamma'_{m})\mathfrak{c'}}, \ldots.
\]
\cite[p.~340]{Dirichlet37}
\end{quote}
This makes it harder to appreciate the nature of the cancellation trick. Moreover, although values $\Theta, \Phi, \Omega, \Omega', \ldots$ can be used to define the individual characters, these tuples and the corresponding representation play no role in the subsequent proof, which depends only on the orthogonality relations and the multiplicative nature of the characters. It seems reasonable, then, to seek a manner of expression that abstracts away the details of the representation. We saw that in his 1841 paper on arithmetic progressions in the quadratic integers, Dirichlet briefly used the expression $S\Omega_n$ to denote the result of summing the values of $\Omega_n$ over all possible combinations of roots that occur in $\Omega$. Kronecker maintained the dependence of the characters on the defining tuples of data, but found a much more elegant notation for expressing the dependence. He denoted the character corresponding to the tuple of parameters $(k)$ by $\Omega^{(k)}$, and in the case of a modulus $m$, he expressed the second orthogonality relation by writing
\[
 \sum_{(k)} \Omega^{(k)}(r_0) = \ph(m),
\]
when $r_0$ is congruent to $1$ modulo $m$, and
\[
 \sum_{(k)} \Omega^{(k)}(r) = 0
\]
otherwise. In his 1883 paper on general characters, Weber adopted a curious means of abstracting the representation of the characters: he simply assigned arbitrary indices to the characters, listing them as $\chi_1, \ldots, \chi_h$. He then expressed the second orthogonality principle without summation notation, as
\[
\chi_{1}(\Theta) + \chi_{2}(\Theta) + \ldots + \chi_{h}(\Theta) = 0,
\]
for each group element $\Theta$. In 1896, however, de la Vall\'ee-Poussin adopted notation $S_\chi$ for summation over characters:
\begin{quote}
Consider \ldots the sum extending over all the characters, that is to say over all the systems of roots
\[
S_{\chi}\chi(n)=S_{\omega}\omega_{1}^{\nu_{1}}\omega_{2}^{\nu_{2}} \ldots
\]
\ldots\ \emph{For every number $n$, the sum extending over the totality of characters satisfies
\[
S_{\chi}\chi(n)=0,
\]
the only exception being the case where
\[
n\equiv 1 \pmod{M},
\]
because then all the indices are zero and one has
\[
S_{\chi}\chi(n)=\ph(M).
\]}
\cite[pp.~14--15]{Poussin96}
\end{quote}
It is notable that he chose a symbol distinct from the usual summation symbol, $\sum$, which he used for sums ranging over natural numbers. Nonetheless, he seems to be the only nineteenth century author to have taken summation over characters at face value. 

Setting aside the orthogonality relation, let us consider the subsequent calculation, involving the $L$-series, where those identities are put to use. We have observed that the modern notation $L(s, \chi)$ allows us to express the dependence of an $L$-series on the character $\chi$, and that the notation $\sum_\chi \overline{\chi(m)}\log L(s, \chi)$ allows us to sum over characters, but these means of expression were not available to Dirichlet. In the case of a prime modulus $p$, Dirichlet defined the $L$ series
\[
 L_0, L_1, \ldots, L_{p-2},
\]
where the index corresponds to a particular numeric parameter occurring in the algebraic expression that we now recognize as the value of the corresponding character, and 
\[
\log L_{0} + \Omega^{-\gamma_{m}} \log L_{1} + \Omega^{-2\gamma_{m}}\log L_{2} + \ldots + \Omega^{-(p-1)\gamma_{m}} \log L_{p-2}
\]
to sum over the $p - 1$ many $L$ series in the case of a prime modulus. In the case of a general modulus $k$, each $L$ series has a similar denotation
\[
L_{\mathfrak{a},\mathfrak{b},\mathfrak{c},\mathfrak{c}', \ldots} 
\]
where $\mathfrak{a},\mathfrak{b},\mathfrak{c},\mathfrak{c}', \ldots$ are a sequence of numeric parameters that appear in the algebraic expression for the general character, and the summation is denoted 
\[
\sum \Theta^{-\alpha_{m}\mathfrak{a}} \Phi^{-\beta_{m}\mathfrak{b}}\Omega^{-\gamma_{m}\mathfrak{c}}\Omega^{-\gamma'_{m}\mathfrak{c}'} \cdots \log L_{\mathfrak{a},\mathfrak{b},\mathfrak{c},\mathfrak{c}', \ldots}
\]
where the summation ranges over the $\ph(k)$ many choices of values of $\mathfrak{a}, \mathfrak{b}, \mathfrak{c}, \mathfrak{c'}, \ldots$. 
Thus Dirichlet took the $L$ series to depend on particular tuples of numeric parameters involved in the definition of the characters, and took summations to range over these parameters. Dedekind's 1863 presentation followed Dirichlet in this respect, as did de la Vall\'ee-Poussin's 1897 presentation. Hadamard in 1896 and Landau in 1909 adopted a tack similar to Weber's, assigning arbitrary indices to the characters, and then letting the $L$-series depend on those indices. For example, Hadamard wrote $\psi_1, \psi_2, \ldots, \psi_{\ph(k)}$ for the list of characters modulo $k$, and defined the $L$-functions as follows:
\[
L_{v}(s)=\sum_{n=1}^{\infty}\frac{\psi_{v}(n)}{n^{s}}.
\]
The key summation over the characters is then written $\sum_{v}\frac{\log L_{v}(s)}{\psi_{v}(m)}$.

To the modern eye, it seems strange to assign otherwise meaningless indices to the characters in order to express the functional dependence of the $L$ series on a character and to sum over them, when one can just write $L(s, \chi)$ and $\sum_\chi$. But while it was perfectly natural in the nineteenth century to sum over integers, summing over the functions themselves may not even have occurred to these authors. It is not until 1897 that we first see $L$ series expressed as a functional dependence on characters, when de la Vall\'ee Poussin introduced the notation $Z(s, \chi)$. Subsequent authors adopted the notation $L(s, \chi)$, reverting back to Dirichlet's use of the letter $L$. By 1927, for example, Landau was using $L(s, \chi)$ and $\sum_\chi$ just as we do today, and from then on the usage seems to have stuck.

\subsection{Extensionalization}

Let $f(x)$ be the function on the real numbers defined by $f(x) = 3 x^2 + 1$. In logical parlance, the \emph{intension} of this last expression is the manner of presentation, in some sense --- if not the purely syntactic string of symbols, something close to it. In contrast, the \emph{extension} is the abstract object denoted, that is, the abstract input-output relation. Today, when we refer to functions, we generally have their extensions in mind. A note of intensionality creeps in when we say things like ``the leading coefficient of $f$'' or ``the constant term of $f$,'' but when called on to explain what we mean, we are generally able to clarify the fact that by ``$f$'' we really mean the expression for $f$ rather than the object itself. The extensional nature of the function concept is embodied in the fact that when we define a functional $F(f)$ on a collection of functions, we ensure the definition does not depend on the manner of presentation of $f$, since $F$ is supposed to ``act'' on the extension, not the intension. 

In ``Concept,'' we argued that this distinction was not as clearly drawn in the nineteenth century treatment of functions. Early instances of functions --- not just functions on the real and complex numbers, but also objects like permutations, automorphisms, and so on --- were more tightly associated with a manner of expression. The history of the treatment of characters in Dirichlet's theorem shows exactly this sort of ambiguity, and a gradual move towards an extensional treatment.\footnote{A referee has suggested that ``abstraction'' and ``abstract treatment'' may be more apt than ``extensionalization'' and ``extensional treatment,'' since ``extensionality'' is often associated with a set-theoretic interpretation of functions. As the referee concedes, however, that the word ``abstract'' has multiple connotations, and so we have stuck with the more focused terminology.} 

Consider, for example, the definition of the concept of character itself. For each $k$, the set of characters modulo $k$ can be defined extensionally, as the set of nonzero homomorphisms from $(\mathbb Z/ k \mathbb Z)^*$ to the complex numbers, or intensionally, as functions defined by certain algebraic expressions involving certain primitive elements modulo the prime powers occurring in the factorization of $k$, and certain complex roots of unity. Even though the two definitions give rise to the same set of characters, proofs can differ in the extent to which they rely on the specific representations or the abstract characterizing property. Dirichlet's proof relied only on the symbolic representations, but we have seen that later proofs emphasized the key properties of the characters, which were extensional in nature. 

Recall also that Dirichlet divided the $L$ series into three classes, depending on a corresponding division of the characters on which they depend. Dirichlet described the division in terms of the tuples of roots appearing in the algebraic expressions, whereas a modern characterization describes the three kinds of characters as follows:
\begin{enumerate}
 \item the character with constant value $1$
 \item the (other) real-valued characters
 \item the (other) complex-valued characters
\end{enumerate}
What is perhaps surprising is that even as later authors introduced notation like $\chi$ or $\psi_i$ to range over characters, they still carried out the classification in terms of the roots. For example, both Dedekind's and Hadamard's division of the characters into the trivial, real, and complex cases was also described in terms of the characters' representations, even though the distinction is naturally expressed in terms of the values they take. Kronecker and de la Vall\'ee-Poussin provided both descriptions, and even though Kronecker made it clear that all operations and classifications can be carried out, algorithmically, in terms of the canonical representations, his careful choice of notation and organization made the extensional properties salient. By 1927, Landau clearly favored the extensional characterization in his textbook.

As yet another means of highlighting the difference between intensional and extensional ways of thinking about functions, we will close this section by noting that a number of the authors we considered adopted a strikingly similar means of describing identities parameterized by the characters. Recall that after stating the generalized version of the Euler product identity (\ref{dirichletEulergeneral}), Dirichlet wrote:
\begin{quote}
The general equation, in which the different roots $\theta, \varphi, \omega, \omega', \ldots$ can be combined with one another arbitrarily, clearly contains $K$-many particular equations.
\end{quote}
The notion of a single identity ``containing'' $K$-many particular equations sounds strange to us today. In contrast to thinking of an identity like $e^{x + y} = e^x + e^y$ as a single equation in which $x$ and $y$ are taken to range over the real or complex numbers, it is almost as though Dirichlet conceived of the generalized Euler product formula as a \emph{template}, or a \emph{schema}, for the particular assertions obtained by instantiating the variables $\theta, \ph, \omega, \omega', \ldots$ with the particular data representing each character. In a similar way, when Dedekind defined the $L$ series in 1863, he wrote:
\begin{quote}
Since these roots can have $a, b, c, c', \ldots$ values, respectively, the form $L$ contains altogether
$abcc'\ldots = \ph(k)$ different particular series\ldots
\end{quote}
This manner of speaking persisted even after authors began using a single symbol $\chi$ to stand for an arbitrary character. For example, in 1882, Weber, after deriving a pair of identities involving an arbitrary character $\chi$, wrote:
\begin{quote}
Each of the formulas \ldots represents $h$ different formulas, corresponding to the $h$ different characters $\chi_1, \chi_2, \ldots, \chi_h$.
\end{quote}
And in a very similar situation, de la Vall\'ee Poussin wrote in 1897:
\begin{quote}
\ldots this equation (E) represents in reality $\ph(M)$ distinct ones, which result from exchanging the characters amongst themselves.
\end{quote}
Such language suggests that, to some extent, authors thought of the act of ``instantiating'' a general identity involving characters at a particular character as somewhat different from instantiating a general identity over numbers at a particular number.

\section{Methodology and ontology revisited}
\label{analysis:section}

Let us review some of the general historical trends we have discerned in the treatment of characters. Over time, authors isolated certain symbolic expressions appearing in Dirichlet's proof, viewed them as functions of an integer parameter (or equivalence class) $n$, and baptized them ``characters.'' They isolated important properties of the characters and articulated them in a way that renders them independent from the rest of the proof. Collaterally, this made it possible to generalize the notion of a character on a multiplicative group of residues to the notion of a character on any abelian group.

Initially, each character was seen to be represented by a bundle of defining data, so what we now characterize as a functional dependence on the character was expressed as a dependence on the bundle of data, and a summation over the characters was expressed as a summation of a range of values of the bundle of data. But, over time, the role that the representing data had to play in the proof was diminished. Authors began to adopt notation and patterns of argumentation that suppressed that information, for example, by assigning arbitrary indices to the characters and letting expressions depend on those indices. Ultimately, authors simply began expressing functional dependences on, and summing over, the characters themselves. 

Avoiding the need to refer to any particular representation of the characters meant relying instead on properties of the characters that can be expressed in terms of the values they take on suitable inputs. In other words, it amounted to adopting an extensional view of the characters, in which statements about the characters are cast purely in those terms. In contemporary proofs of Dirichlet's theorem, this is taken to the extreme when we define the set of characters as the set of nonzero homomorphisms from the group in question to the complex numbers, and carry out the proof without indicating any way of representing individual characters, let alone means of computing with them. 

One might describe these changes as ``merely notational,'' or ``merely pragmatic.'' But dismissing them in that way belies the fact that these changes reflect a fundamentally different way of talking about, thinking about, and reasoning about the characters. And this was by no means an isolated example. As we have noted in the introduction, during the nineteenth century the treatment of other mathematical entities that we now take to be instances of sets, functions, or structures evolved in similar ways, and for similar reasons. So the history we have traced here is but one instance of a general transformation in mathematical thought, with a new conception of the basic objects of mathematics and appropriate means of reasoning about them. It seems strange to resist seeing this as a change in ontology. (Gray \cite{gray:92} nicely emphasizes this point.)

According to the historical model described in Section~\ref{metaphysics:section}, we should view the history of Dirichlet's theorem as a response to fundamental methodological pressures, as mathematicians struggled to meet both intrinsic and extrinsic mathematical goals while respecting intrinsic and extrinsic methodological constraints. As philosophers, we should not be interested so much in the historical and psychological contingencies that shaped the process, but, rather, the sense in which the outcome is rational and justified. In other words, we wish to understand the extent to which the methods of contemporary mathematics serve to achieve our mathematical goals, given some conception of those goals and what it means to do mathematics. Attention to the history can bring some of the goals and constraints to light, but then we are left to weigh their importance and assess the merits of the present solution. This is the point at which philosophical analysis must come into play. 

In broad terms, here we will view mathematics as a process by which finite beings attempt to impose a useful order on the complex and varied data that confronts them. The philosophical task is then to develop more refined characterizations of the mathematical process, in terms that adequately reflect the constraints we face as mathematical agents and the goals we pursue. In ``Concept,'' we provided a detailed discussion of some of the various methodological benefits and concerns that accrue to the use of the modern function concept. Let us briefly review these here, and see what they have to tell us about the nature of mathematics.

Treating characters as objects, in all the senses described in Section~\ref{functions:section}, brings a number of methodological benefits. Expressions become simplified, meaning that the reader has to keep track of less information when parsing them, and the author of a proof can record and convey the relevant information more compactly. Proofs become simplified as well, meaning that readers have to keep track of less information while following the argumentative structure of a proof, and authors have to keep track of less information while working out the details. Information that is irrelevant to the argument at hand, or can be made so, is suppressed, making key data and relationships more salient. 

Moreover, proofs became more modular, as properties of the characters were abstracted away and proved separately. This further supports the aim of reducing the amount of information in play at any given point. While developing a theory of the characters, we need only work with their defining properties, and when checking that particular instances of functions are characters, we need only check that these instances satisfy the defining properties. Then, when reasoning about these particular characters, we can invoke results from the general theory, such as the orthogonality lemma, as ``black boxes.'' The fact that extraneous information has been filtered out means that expressions depend on fewer parameters, and inferences depend on fewer assumptions. This makes it easier to check details and avoid mistakes.

Modularity brings additional benefits, in that definitions and theorems that have been abstracted away from the body of the proof can be reused elsewhere. The process of abstraction clarifies the data that serves to parametrize a definition and the hypotheses that are required to establish a proposition. This facilitates not only using the definitions and proposition in other contexts, but also  modifying the definitions and propositions by varying the parameters and hypotheses accordingly. In this way, modularity supports generality as well as reuse. 

Thus, with a modular structuring, dependencies between mathematical components are minimized, and the mathematics becomes easier to understand. It also becomes easier to ensure correctness, and components can be modified and reused. Notice, incidentally, that these are exactly the benefits associated with modularity in software engineering.\footnote{This is a topic is explored in greater detail in \cite{avigad:unp}.}

The key point is that treating characters as objects supports this modularity. To start with, identifying characters as ``things'' means that they can be objects of study. We can make assertions about them, and specify predicates and functions that take them as arguments. Moreover, notations, definitions, and theory designed to handle other ``things'' now applies: we can form sums that range over the characters and reason about them; we can form sets and sequences of characters and reason about them; we can consider groups of characters and reason about them; and so on. In short, all of methods that are available to us for reasoning about mathematical objects become applicable to reasoning about characters.

Given the apparent benefits of treating characters as full-blooded objects, why did it take so long for the mathematical community to do so? When we look back at the history of mathematics, it is hard to appreciate the difficulties that accompany significant shifts in method, but they are substantial. Mathematics is a communal activity: when a mathematician writes a proof, his or her intention is that others will read it and judge it to be informative and correct. This requires that the author and the reader have a common understanding not only as to what is permissible, but also as to what is appropriate and desirable.

In Section~\ref{metaphysics:section}, we enumerated some of the concerns that arise when new methods are introduced. In ``Concept,'' we explored the way these concerns apply specifically to the modern treatment of characters, and to functions more generally. To start with, it is important that the new manner of speaking about functions come with clear rules of use. If there is no agreement as to which inferences are permissible --- for example, under what conditions it is legitimate to consider two expressions denoting functions as ``equal,'' and to substitute one expression for another in a given context --- then the mathematical enterprise falls apart, and mathematicians cannot read each others' proofs.

Moreover, whether the rules of use are presented explicitly or implicitly, there is also the question as to whether they are consistent.  Even if we think of the new treatment of characters as a mere short cut to establishing Dirichlet's theorem, such short cuts are clearly illegitimate if they lead to false or nonsensical conclusions. It is by no means apparent that there are no hidden pitfalls in quantifying over characters, summing over characters, and treating characters as arguments to other functions. It would be mathematically reckless to adopt these devices out of sheer convenience, without some assurances that the results obtained are reliable. As suggested in Section~\ref{metaphysics:section}, to some extent it helps to know how the new methods can be interpreted in terms of the prior methods, bolstering the understanding that \emph{if} we try to view talk of characters as short cuts to proving new theorems, the long way is still, in principle, open to us.

Even if the new rules of use seem to be reliable, there is still the question as to whether they are meaningful. We argued in Section~\ref{transition:section} that early authors tended to think of characters as symbolic expressions of a certain kind, or at least, as entities with canonical representations as such symbolic expressions. If the new methods no longer support such a view, one has to come to terms with the question of how one \emph{should} think of a character. Put succinctly, once we have proved a statement about characters, what do we know?

And even if we come to believe that a certain manner of working with characters is consistent, legitimate, and meaningful, there is still the question as to whether it constitutes \emph{good mathematics}, which is to say, whether it furthers our epistemic goals and provides satisfactory answers to our questions. This issue becomes pressing when we try to reconcile a computational conception of mathematics with the new methods of abstraction. For most of its history, mathematics was essentially computational, supplying methods of calculation that could be used to predict the motion of the planets, succeed in games of chance, and compute lengths and magnitudes of all sorts. A central feature of the modern treatment of characters is that it suppresses details of how to represent and compute with individual characters, and often even eliminates these details entirely. We may feel as though we have an understanding of what it means for a function, viewed as a general procedure, to take a natural number as input, but what does it mean for a function to take a character, viewed abstractly, as input? If we expect a mathematical theory of characters to tell us how to represent them and compute with them, then a theory that fails to provide that information is simply defective.

Separating concerns as we have done here is somewhat artificial. For example, maintaining a computational view of characters is one way of interpreting their meaning, and the ability to ascribe any sort of meaning to mathematical objects tends to clarify the rules of use and support the belief in these rules are consistent. Notice, also, that on our analysis, the factors that ultimately support adopting a modern treatment of functions are an uneasy mix of pragmatic, empirical, and broadly philosophical considerations. That does not mean that they are not good reasons, however, nor that we have not made important philosophical progress by understanding them better.

In the latter half of the nineteenth century, Frege's development of formal logic was designed to represent mathematical language and methods of reasoning, and offer clear recommendations as to proper usage. Famously, the notion of ``function'' is central to his account, as well as an understanding of the relationship between ``function'' and ``object.'' In the remaining sections of this essay, we will consider Frege's analysis, and argue that his logical and philosophical choices were influenced by many of the same considerations that were faced by his mathematical peers.

\section{Frege's view of functions and objects}
\label{frege:section}

In 1940, Alonzo Church presented a formulation of type theory \cite{church:40}, now known as ``simple type theory.'' Simple type theory can serve as a foundation for a significant portion of mathematics, and, indeed, is the axiomatic foundation of choice for a number of computational interactive theorem provers today \cite{gordon:melham:93,harrison:07c,nipkow:et:al:02}. One starts with some basic types, say, a type $\BB$ of Boolean truth values and a type $\NN$ of natural numbers, and one forms more complex types $\sigma \times \tau$ and $\sigma \to \tau$ from any two types $\sigma$ and $\tau$. Intuitively, elements of type $\sigma \times \tau$ are ordered pairs, consisting of an element of type $\sigma$ and an element of type $\tau$, and elements of type $\sigma \to \tau$ are functions from $\sigma$ to $\tau$. In a type-theoretic approach to the foundations of mathematics, one identifies sets of natural numbers with predicates, which is to say, elements of type $\NN \to \BB$. Binary relations on the natural numbers are then elements of type $\NN \times \NN \to \BB$, and sequences of natural numbers are elements of type $\NN \to \NN$. Objects at this level are called \emph{type 1} elements, because they require one essential use of the function space arrow. Integers can be identified as pairs of natural numbers and rationals can be identified as pairs of integers in the usual ways. Real numbers are then Cauchy sequences of rationals (elements of type 1), or equivalence classes of such, which puts them at type 2. Functions from the reals to the reals and sets of reals are then elements of type 3, and sets of functions from the reals to reals or collections of sets of real numbers are then elements of type 4. For example, the collection of Borel sets of real numbers is an element of type 4, as is Lebesgue measure, which maps certain sets of real numbers to the real numbers. A set of measures on the Borel sets of the real numbers is an element of type 5. And so on up the hierarchy.

Simple type theory can be viewed as a descendant of the ramified type theory of Russell and Whitehead's \emph{Principia Mathematica} \cite{russell:whitehead:10}, which, in turn, was inspired by the formal system of Frege's \emph{Grungesetze der Arithmetik} \cite{frege:grundgesetze}. Starting with a basic type of individuals, Frege's system also has variables ranging over higher-type functionals, and so can be seen as an incipient form of modern type theory. For that reason, it may come as a surprise to logicians familiar with the modern type-theoretic understanding that the foundational outlook just described is \emph{not at all} the image of mathematics that Frege had in mind. It is this image that we wish to explore here.

Frege took concepts to be instances of functions; for example, in ``Function and concept'' he wrote that ``a concept is a function whose value is always a truth value'' \cite[p.~139]{Frege91}.\footnote{We should note that in this section we will focus on his views from 1884 onwards.  Prior to this, he seems to have held a different view of concepts, though he still maintained that they are not objects; see \cite[p.~136]{compthought}.} And, throughout his career, he was insistent that functions are not objects. The third ``fundamental principle'' in his \emph{Grundlagen der Arithmetik} of 1884 was ``never to lose sight of the distinction between concept and object''
\footnote{``\ldots der Unterschied zwischen Begriff und Gegenstand ist in Auge zu behalten.''} 
\cite[Introduction]{Frege84}, and he later asserted that ``it will not do to call a general concept word the name of a thing'' \cite[\S\,51]{Frege84}.\footnote{``\ldots ist es unpassend, ein allgemeines Begriffswort Namen eines Dinges zu nennen.''} The distinction features prominently in his essays ``Function and concept,'' ``Comments on \emph{Sinn} and \emph{Bedeutung}'' and ``Concept and object'' of 1891, 1891/2, and 1892, respectively.

According to Frege, the proper distinction is tracked by linguistic usage: objects are denoted by words and phrases that can fill the subject role in a grammatical sentence, whereas concepts are denoted by words and phrases that can play the role of a predicate. In ``Concept and object'' he wrote:
\begin{quote}
 We may say in brief, taking ``subject'' and ``predicate'' in the linguistic sense: a concept is the \emph{Bedeutung} of a predicate; an object is something that can never be the whole \emph{Bedeutung} of a predicate, but can be the \emph{Bedeutung} of a subject.\footnote{``Wir k\"{o}nnen kurz sagen, indem wir ``Pr\"{a}dikat'' und ``Subjekt'' im sprachlichen Sinne verstehen: Begriff ist Bedeutung eines Pr\"{a}dikates, Gegenstand ist, was nie die ganze Bedeutung Pr\"{a}dikates, wohl aber Bedeutung eines Subjekts sein kann.'' The word \emph{Bedeutung} is often translated as ``reference'' or ``denotation.'' But for difficulties in the translation, see \S4 of the introduction to Beaney \cite{FregeReader}.} 
 \cite[pp.~198]{cando}
\end{quote}
And:
\begin{quote}
 A concept---as I understand the word---is predicative. On the other hand, a name of an object, a proper name, is quite incapable of being used as a grammatical predicate.\footnote{``Der Begriff---wie ich das Wort verstehe---ist pr\"{a}dikativ.  Ein Gegenstandsname hingegen, ein Eigenname ist durchaus unf\"{a}hig, als grammatisches Pr\"{a}dikat gebraucht zu werden.''} 
 \cite[pp.~193]{cando}
\end{quote}
In the sentence, ``Frege is a philosopher,'' the word ``Frege'' denotes an object, and the phrase ``is a philosopher'' denotes a concept. Frege clarified the distinction by explaining that functional expressions, including concept expressions, are ``unsaturated,'' or incomplete. These stand in contrast to signs that are used to denote objects, which are complete in and of themselves. For example, in the sentence ``Frege is a philosopher,'' the expression ``Frege'' is saturated, and succeeds in picking out an object. In contrast, the expression ``\ldots is a philosopher'' contains a gap, and fails to name an object until one fills in the ellipsis, at which point the expression denotes a truth value.\footnote{While the distinction between saturated and unsaturated expressions is cast as a distinction between linguistic signs, in his 1904 essay ``What is a Function?'' Frege made it clear that the dichotomy extends to functions and objects themselves: ``The peculiarity of functional signs, which we here called `unsaturatedness', naturally has something answering to it in the functions themselves. They too may be called `unsaturated' \ldots'' 
(``Der Eigent\"{u}mlichkeit der Fuktionszeichen, die wir Unges\"{a}ttigtheit genannt haben, entspricht nat\"{u}rlich etwas an den Funktionen selbst.  Auch diese k\"{o}nnen wir unges\"{a}ttigt nennen \ldots'') 
\cite[p.~665]{Frege04}.}

Having distinguished between concepts and objects in such a way, Frege had to deal with objections, such as the one he attributed to Benno Kerry in ``Concept and object.'' In the sentence ``The concept `horse' is a concept easily attained'' the concept denoted by ``horse'' does fill the subject role. Frege's surprising answer was to deny that the phrase ``the concept `horse' '' denotes a concept. He conceded that this sounds strange:
\begin{quote}
 It must indeed be recognized that we are confronted by an awkwardness of language\ldots if we say that the concept \emph{horse} is not a concept\ldots.\footnote{``Es kann ja nicht verkannt werden, da{\ss} hier eine freilich unvermeidbare sprachliche H\"{a}rte vorliegt, wenn wir behaupten: der Begriff Pferd ist kein Begriff \ldots.''} 
 \cite[pp.~196--197]{cando}
\end{quote}
Yet, he insisted, this is what we must do. He was already clear about this in the \emph{Grundlagen}:
\begin{quote}
 The business of a general concept word is precisely to signify a concept.  Only when conjoined with the definite article or a demonstrative pronoun can it be counted as the proper name of a thing, but in that case it ceases to count as a concept word.  The name of a thing is a proper name.\footnote{``Ein allgemeines Begriffswort bezeichnet eben einen Begriff.  Nur mit dem bestimmten Artikel oder einem Demonstrativpronomen gilt es als Eigenname eines Dinges, h\"{o}rt aber damit auf, als Begriffswort zu gelten.  Der Name eines Dinges ist ein Eigenname.''} 
 \cite[\S 51]{Frege84}
\end{quote}
And so, in ``Concept and object,'' he reminded us:
\begin{quote}
 If we keep it in mind that in my way of speaking expressions like ``the concept $F$'' designate not concepts but objects, most of Kerry's objections already collapse.\footnote{``Wenn wir festhalten, da{\ss} in meiner Redeweise Ausdr\"{u}cke wie ``der Begriff $F$'' nicht Begriffe, sondern Gegenst\"{a}nde bezeichnen, so werden die Einwendungen \emph{Kerrys} schon gr\"{o}{\ss}tenteils hinf\"{a}llig.''} 
 \cite[pp.~198--199]{cando}
\end{quote}
He similarly urged us to reconstrue expressions like ``all mammals have red blood'' as ``whatever is a mammal has red blood'' so as to avoid the impression that the predicate ``has red blood'' is being applied to an object, ``mammal.'' Although these examples deal with concepts, Frege's analysis makes it clear that he intended the linguistic separation to remain operant for other kinds of functions as well.

At the same time, Frege was equally dogmatic in insisting that what we commonly take to be mathematical objects really \emph{are} mathematical objects as such. The introduction to his \emph{Grundlagen} begins as follows:
\begin{quote}
 When we ask someone what the number one is, or what the symbol 1 means, we get as a rule the answer ``Why, a thing.''\footnote{``Auf die Frage, was die Zahl Eins sei, oder was das Zeichen 1 bedeute, wird man meistens die Antwort erhalten: nun, ein Ding.'' All our
 translations from the \emph{Grundlagen} are taken from the Austin translation cited in the references.} \cite[Introduction]{Frege84}
\end{quote}
The claim is so curious as to give one pause.\footnote{We are grateful to Steve Awodey for this observation.} The fact that Frege used such a brazen rhetorical flourish to frame the whole project makes it clear just how central the issue is to his analysis. Once again, he took the distinction to be tracked by linguistic use. For example, because the number $7$ plays the role of a subject in the statement ``$7$ is odd,'' $7$ must be an object. But, once again, Frege had to deal with sentences where the syntactic role of a number is murkier. For example, he considered uses of number terms in language that are attributive and do not occur prefixed by the definite article, for example, ``Jupiter has four moons'' \cite[\S 57]{Frege84}.  He wrote
\begin{quote}
 ``\ldots our concern here is to arrive at a concept of number usable for the purposes of science; we should not, therefore, be deterred by the fact that in the language of everyday life number appears also in attributive constructions.  That can always be got round.''\footnote{``Da es uns hier darauf ankommt, den Zahlbegriff so zu fassen, wie er f\"{u}r die Wissenschaft brauchbar ist, so darf es uns nicht st\"{o}ren, dass im Sprachgebrauche des Lebens die Zahl auch attributiv erscheint.  Das l\"{a}sst sich immer vermeiden.''} \cite[\S 57]{Frege84}
\end{quote}
Specifically, it can be got round by writing an attributive statement such as ``Jupiter has four moons'' as ``the number of Jupiter's moons is the number 4, or 4'' \cite[\S 57]{Frege84}, thereby eliminating the attributive usage.

So, for Frege, functions are not objects, but numbers are, because they play the subject role in mathematical statements and can be used with the definite article. There is clearly a difficulty lurking nearby. At least from a modern standpoint, we tend to view functions, sequences, sets, and structures as objects, and certainly in Frege's time locutions such as ``the function $f$'' and ``the series $s$'' were common. Frege's response was similar to his response to Kerry's objection, namely, to deny that that expressions like these denote functions.
To understand how this works, consider the fact that Frege's logical system includes an operator which takes any function $f$ from objects to objects and returns an object, $\overset{,}{\varepsilon} f (\varepsilon)$, intended to denote its ``course-of-values'' or ``value range.'' If $f$ is a concept, which is to say, a function which for each object return a truth value, the course-of-values of $f$ is called the ``extension'' of the concept. Frege's Basic Law V asserts that two functions which are extensionally equal---that is, which return equal output values for every input---have the same courses-of-values.

Frege used these courses-of-values and extensions as object-proxies for functions and concepts. This is how he analyzed the concept of a cardinal number. Let $F$, for example, be a second-level concept, such that $F$ holds of a first-level concept $f$ if and only $f$ holds of exactly one object. Frege took the number one to be the extension of $F$, thereby achieving the goal of making the number one, well, a thing. But this ``pushing down'' trick is central to the methodology of the \emph{Grundgesetze}: whenever the formal analysis of common mathematical objects seems to suggest identifying such objects as functions or concepts, Frege avoided doing so by replacing the function or concept with its extension. For example, in the \emph{Grundgesetze} he circumvented the need to define mathematical operations on sequences and relations  construed as functions, defining the operations rather on the associated courses-of-values.\footnote{In fact, the definition of the number one earlier in the paragraph describes, more precisely, the construction in the \emph{Grundlagen}. In the \emph{Grundgesetze}, he took $F$ to be a \emph{first}-level concept that holds of classes, i.e.~extensions of concepts, that contain one element. This is a nice illustration of how the ``pushing down'' trick can be used repeatedly to avoid the use of higher types. See Reck \cite[Section 5]{reck:07} for a discussion of the two definitions, and Burgess \cite{burgess:05} for an overview of Frege's methodology.\label{number:footnote}}  Referring to Frege's concepts as ``attributes'' and their extensions as ``classes,'' Quine described the difference as follows:
\begin{quote}
Frege treated of attributes of classes without looking upon such discourse as somehow reducible to a more fundamental form treating of attributes of attributes. Thus, whereas he spoke of attributes of attributes as \emph{second-level} attributes, he rated the attributes of classes as of first level; for he took all classes as rock-bottom objects on par with individuals. \cite[p.~147]{quine:55}
\end{quote}

Frege never got so far as developing mathematical analysis in his system, and we cannot say with certainty how he would have developed, for example, ordinary calculus on the real numbers. But there is a strong hint that here, too, he would have taken, for example, operations like integration and differentiation to operate on extensions, rather than functions, in his system. He touched on the history of analysis in his \emph{Function and concept} of 1891, and noted that, for example, differentiation can be understood as a higher-type functionals.
\begin{quote}
 Now at this point people had particular second-level functions, but lacked the conception of what we have called second-level functions. By forming that, we make the next step forwards. One might think that this would go on. But probably this last step is not so rich in consequences as the earlier ones; for instead of second-level functions one can deal, in further advances, with first-level functions---as shall be shown elsewhere.\footnote{``Damit hatte man nun einzelne Funktionen zweiter Stufe, ohne jedoch das zu erfassen, was wir Funktion zweiter Stufe genannte haben.  Indem man dies tut, macht man den n\"{a}chsten Fortschritt.  Man k\"{o}nnte denken, dass dies so weiter ginge.  Wahrscheinlich ist aber schon dieser letzte Schritt nicht so folgenreich wie die fr\"{u}heren, weil man statt der Funktionen zweiter Stufe im weiteren Fortgang Funktionen erster Stufe betrachten kann, wie an einem anderen Orte gezeigt werden soll.''} \cite[p.~31]{Frege91}
\end{quote}
Presumably, he had the method of replacing functions by their extensions in mind.

Notice, incidentally, that Frege's method of representing mathematical functions as courses-of-values has the effect that mathematical functions are treated extensionally. For example, defining the integral as an operation that applies to a course-of-values means that integration cannot distinguish mathematical functions that are extensionally equal, since any two descriptions of a function that satisfy extensional equality have the same course-of-values, by Basic Law V.

There are other interesting features of Frege's treatment of functions that push us away from identifying them with the functions of ordinary mathematics. For example, for Frege, every function has to be defined on the entire domain of individuals; even if one is interested in the exponential function on the real numbers, one has to specify a particular (but arbitrary) value of this function for every object in existence.\footnote{However, Patricia Blanchette has argued \cite{blanchette:12} that Frege intended theories presented in his formal system to treat objects in the domain of a particular subject, in which case ``every object in existence'' really means ``every object in the theory's intended domain.''} And the separation of functions and objects has other effects on the system. There is only one basic type, so, for example, truth values live alongside everything else. There is no notion of identity between higher-type objects---the equality symbol can only be applied to equality between objects---even though Frege pointed out that one can define an extensional notion of ``sameness'' of functions and concept, for example, saying two functions from individuals to individuals are the ``same'' if their values are identical at each input. Frege's system, of course, includes the axiom of universal instantiation. In contemporary notation, this would be expressed as $\fa x \ph(x) \limplies \ph(a)$ where $x$ is a variable ranging over individuals and $a$ is any individual term. It also includes the corresponding axiom $\fa F \ph(F) \limplies \ph(A)$, where $F$ ranges over functions from objects to objects. Notably, however, the system does not include analogous axioms for elements of the higher types: the ``pushing down trick'' obviates the need for these.

All things considered, Frege's foundational treatment of mathematics seems closer to modern set-theoretic treatments, where there is one homogeneous universe of individuals. Truth values are individuals, numbers are individuals, mathematical sequences and series are individuals---all bona-fide mathematical objects are individuals. Functions are special sorts of entities that our partial expressions refer to when we make statements about objects, but they are not objects in their own right. As Marco Panza puts it:
\begin{quote}
 \ldots according to Frege, appealing to functions is indispensable in order to fix the way his formal language is to run, but functions are not as such actual components of the language. More generally, functions manifest themselves in our referring to objects---either concrete or abstract---and making statements about them, but they are not as such actual inhabitants of some world of \emph{concreta} and \emph{abstracta}. \cite[p.~14]{panza:unp}
\end{quote}
This is not to say that functions are any less ``real'' or objective than mathematical objects like numbers, only that they play a distinct role: they allow us to define objects, say things about objects, and reason about objects, but they are not objects themselves.

\section{Frege's foundational concerns}
\label{frege:section:b}

We have seen that a curious tension lies at the core of Frege's formal representation of mathematics. On the one hand, Frege asserted, repeatedly, that functions, in the logical and linguistic sense, are not objects. On the other hand, when it comes to formalizing mathematical constructions, he clearly felt that functions, in the mathematical sense, \emph{have to be} objects. His course-of-values operator, together with his Basic Law V, allowed him to have his cake and eat it too, maintaining clear borders between the two realms while passing between them freely. But Frege is often taken to task for failing to realize that this strategy opens the door to Russell's paradox. Indeed, the strategy feels like a hack, a desperate attempt to satisfy the two central constraints. Why was he so committed to them? The goal of this section is to suggest that the concerns Frege was trying to address with the design of the logic of the \emph{Grundgesetze} parallel some of the informal mathematical concerns we were able to discern in the nineteenth century treatment of characters.

When one speculates as to the philosophical and logical considerations that influenced the design of Frege's logic, two possibilities come to mind. One is that Frege determined that functions and objects should be separate on broad ontological grounds, and then designed the logic accordingly. The other is that he designed the logic, determined it worked out best with a separation of individuals and functions, and read off the ontological stance from that. But, in fact, there is no clear distinction between these two descriptions. Frege designed his logic to try to model scientific practice at its best, and account for and support its successes while combating and eliminating confusions. The examples in the previous section show that Frege had no qualms about reinterpreting ordinary locutions and reconstruing everyday language, so he was by no means slave to naive ontological intuitions. But even when doing so he appealed to intuitions to convince us that the reconstruals are reasonable. Thus ``doing ontology'' meant analyzing the practice, sorting out intuitions, and trying to regiment and codify them in a coherent and effective way. From the other direction, ``getting the logic to work'' meant being able to account for the informal practice effectively and efficiently, and supporting our intuitions to the extent that they can be fashioned into a coherent system. So it is not a question as to whether the ontology or the logic comes first; working out the ontology and designing the logic are part and parcel of the same enterprise. The following questions therefore seem more appropriate:
\begin{enumerate}
 \item What considerations pushed Frege to maintain the sharp distinction between function and object?
 \item What considerations pushed Frege to identify mathematical entities, including ordinary mathematical functions, as objects?
\end{enumerate}
Let us consider each in turn.

It seems to us that the answer to the first question is simply that Frege felt that failure to respect the distinction results in linguistic confusion.
\begin{quote}
 If it were correct to take ``one man'' in the same way as ``wise man,'' we should be able to use ``one'' also as a grammatical predicate, and to be able to say ``Solon was one'' just as much as ``Solon was wise.'' It is true that ``Solon was one'' can actually occur, but not in a way to make it intelligible on its own in isolation. It may, for example, mean ``Solon was a wise man,'' if ``wise man'' can be supplied from the context. In isolation, however, it seems that ``one'' cannot be a predicate. This is even clearer if we take the plural. Whereas we can combine ``Solon was wise'' and ``Thales was wise'' into ``Solon and Thales were wise,'' we cannot say ``Solon and Thales were one.'' But it is hard to see why this should be impossible, if ``one'' were a property both of Solon and of Thales in the same way that ``wise'' is.\footnote{``Wenn `Ein Mensch' \"{a}hnlich wie `weiser Mensch' aufzufasen w\"{a}re,  so sollte man denken, dass, `Ein' auch als Praedicat gebraucht werden k\"{o}nnte, sodass man wie `Solon war weise' auch sagen k\"{o}nnte `Solon war Ein' oder `Solon war Einer'.  Wenn nun der letzte Ausdruck auch vorkommen kann, so ist er doch d\"{u}r sich allein nicht verst\"{a}ndlich. Er kann z.B. heissen: Solon war ein Weiser, wenn `Weiser' aus dem Zusammenhange zu erg\"{a}nzen ist.  Aber allein scheint `Ein' nicht Praedicat sein zu k\"{o}nnen.  Noch deutlicher zeigt sich dies beim Plural.  W\"{a}hrend man `Solon war weise' und `Thales war weise' zusammenziehen kann in `Solon und Thales waren weise,' kann man nicht sagen `Solen und Thales waren Ein'.  Hiervon w\"{a}re die Unm\"{o}glichkeit nicht einzusehen, wenn `Ein' sowie `weise' eine Eigenschaft sowohl des Solon als auch des Thales w\"{a}re''.} \cite[\S 29]{Frege84}
\end{quote}
In other words, even though in some contexts an object word like ``one'' can \emph{appear} to be used as a predicate, and in other contexts a concept can \emph{appear} to be used as a subject, closer inspection shows that these uses do not conform to the rules that govern the use of prototypical subjects and predicates, and so should not be categorized in the naive way.

One of Frege's favorite pastimes was to show that assertions made by philosophical and mathematical colleagues degenerate into utter nonsense when they fail to maintain sufficient linguistic hygiene. For example, in his  1904 essay, ``What is a Function,'' Frege was critical of conventional mathematical accounts of variables and functions. It is a mistake, he said, to think of a variable as being an object that varies:
\begin{quote}
  \ldots a number does \emph{not} vary; for we have nothing of which we could predicate the variation.  A cube never turns into a prime number; an irrational number never becomes rational.\footnote{``Folglich ver\"{a}ndert sich die Zahl gar nicht; denn wir haben nichts, von dem wir die Ver\"{a}nderung aussagen k\"{o}nnten.  Eine Kubikzahl wird nie zu einer Primzahl, und eine Irrationalzahl wird nie rational.''} \cite[p.~658]{Frege04}
\end{quote}
 He took the mathematician Emanuel Czuber to task for giving such a sloppy account of variables and functions in an introductory mathematical text. For example, he criticized Czuber's terminology ``a variable assumes a number'' \cite[288]{Frege04} as being incomprehensible. On Czuber's account, a variable is an ``indefinite number,'' so the terminology can be rephrased ``an indefinite number assumes a (definite) number''; but where we may talk about an object assuming a property, what can it mean for an object to assume another object?
\begin{quote}
 In other connections, indeed, we say that an object assumes a property, here the number must play both parts; as an object it is called a variable or a variable magnitude, and as a property it is called a value.  That is why people prefer the word ``magnitude'' to the word ``number''; they have to deceive themselves about the fact that the variable magnitude and the value it is said to assume are essentially the same thing, that in this case we have \emph{not} got an object assuming different properties in succession, and that therefore there can be no question of a variation.\footnote{``Sonst sagt man wohl, da{\ss} ein Gegenstand eine Eigenschaft annehme; heir mu{\ss} die Zahl beide Rollen spielen; als Gegenstand wird sie Variable oder ver\"{a}nderliche Gr\"{o}{\ss}e, als Eigenschaft wird sie Wert genannt.  Darum also zieht man das Wort `Gr\"{o}{\ss}e' dem Worte `Zahl'ert, den sie angeblich annimmt, im Grunde dasselbe sind, da{\ss} man gar nicht den Fall hat, wo ein Gegenstand nacheinander verschiedene Eigenschaften annimmt, da{\ss} also von Ver\"{a}nderung in keiner Weise die Rede sein kann.''} \cite[p.~660--661]{Frege04}
 \end{quote}

The essay closes with the following assessment:
\begin{quote}
 The endeavor to be brief has introduced many inexact expressions into mathematical language, and these have reacted by obscuring thought and producing faulty definitions.  Mathematics ought properly to be a model of logical clarity.  In actual fact there are perhaps no scientific works where you will find more wrong expressions, and consequently wrong thoughts, than in mathematical ones.  Logical correctness should never be sacrificed to brevity of expression.  It is therefore highly important to devise a mathematical language that combines the most rigorous accuracy with the greatest possible brevity.  To this end a symbolic language would be best adapted, by means of which we could directly express thoughts in written or printed symbols without the intervention of spoken language.\footnote{``Das Streben nach K\"{u}rze hat viele ungenaue Ausdr\"{u}cke in die mathematische Sprache eingef\"{u}hrt, und diese haben r\"{u}ckwirkend die Gedanken getr\"{u}bt und fehlerhafte Definitionen zuwege gebracht.  Die Mathematik sollte eigentlich ein Muster von logischer Klarheit sein.  In Wirklichkeit wird man vielleicht in den Schriften keiner Wissenschaft mehr schiefe Ausdr\"{u}cke und infolgedessen mehr schiefe Gedanken finden als in den mathematischen.  Niemals sollte man die logische Richtigkeit der K\"{u}rze des Ausdrucks opfern.  Deshalb ist es von gro{\ss}er Wichtigkeit, eine mathematische Sprache zu schaffen, die mit strengster Genauigkeit m\"{o}glichste K\"{u}rze verbindet. Dazu wird wohl am besten eine Begriffsschrift geeignet sein, ein Ganzes von Regeln, nach denen man durch geschriebene oder gedruckte Zeichen ohne Vermittlung des Lautes unmittelbar Gedanken auszudr\"{u}cken vermag.''} \cite[p.~665]{Frege04}
\end{quote}

Frege aimed to give a clear account of the rules that govern proper logical reasoning. Although, in ordinary language, the line between concepts and objects is sometimes blurry, failure to diagnose and manage the blurriness opens the door to nonsensical reasoning. Even though words like ``one'' and ``horse'' sometimes seem to denote both concepts and objects, conflating the two causes problems. For Frege, the only viable solution was to analyze and regiment such uses in a way that cordons off problematic instances. He found that the best way to do this is to maintain a clear separation of concept and object, and then supplement the analysis with an explanation as to how some words seem to cross the divide in certain contexts.

Now let us turn to the second question: why was Frege so dogged in his insistence that mathematical entities like numbers have to be treated as objects, and so persistent, in practice, in pushing mathematical constructions down to that realm? We believe that the answer lies in an observation that we found in Heck \cite{heck:97}: Frege wanted his numbers to be able to count all sorts of entities, and the only way he could make that work was by treating all these entities as inhabitants of the same type. Consider the following statements:
\begin{itemize}
 \item There are two truth values.
 \item There are two natural numbers strictly between 5 and 8.
 \item There are two constant functions taking values among the truth values.
 \item There are two characters on $(\ZZ / 4\ZZ)^*$.
 \item There are two subsets of a singleton set.
\end{itemize}
Frege would have insisted that the word ``two'' in each of these statements refers to the same object. We would like to say that the number of truth values is equal to the number of natural numbers between 5 and 8, but if truth values and numbers were different types of entities, his analysis of number would not do that: even if Frege had a notion of identity for each type, we would have to define a different notion of two for each such type. In other words, for each type $\sigma$, we would have to define a concept $\mathrm{Two}_\sigma$ that holds of concepts of arguments of type $\sigma$ under which two elements fall.\footnote{In Frege's system, in which the equality symbol can only be used with objects, this would have to be expressed instead in terms of a ``sameness'' relation for elements of type $\sigma$, for any $\sigma$ other than the type of objects.} Taking extensions according to Frege's construction would yield an object $2_{\sigma}$ for each $\sigma$. But this results in a proliferation of twos, and since $2_\sigma$ and $2_\tau$ are not guaranteed to be the same object, one would have to exercise great care when reasoning about the relationship between them. This is clearly unworkable.

Instead, Frege designed his numbers to count objects, \emph{simpliciter}: $2$ is the extension of the concept of being a concept of \emph{objects} under which exactly two objects fall.\footnote{See footnote~\ref{number:footnote}.} But this means that if you want to count a collection of things, those things have to be elements of the type of objects. This, in turn, provides a strong motivation to locate mathematical entities of all kinds among the type of objects.

We have seen in Sections \ref{modern:formulation:section} to \ref{analysis:section} that what holds true of counting holds true of other mathematical operations, relations, and constructions as well. Contemporary proofs of Dirichlet's theorem have us sum over finite sets of characters just as we sum over finite sets of numbers. We view the general operation here as summation over a finite set of objects, viewing both characters and numbers as such. Contemporary proofs also have us consider groups of characters, just as we consider groups of residues. Once again, we consider these as instances of the general group concept, with the understanding that a group's underlying set can be any set of objects. This allows us to speak of a homomorphism between any two groups, without requiring a different notion of ``homomorphism'' depending on the type of objects of the groups' carriers.

Characters were not the only mathematical entities studied in the latter half of the nineteenth century that encouraged set-theoretic reification. Gauss' genera of quadratic forms, discussed briefly in ``Concept,'' also bear a group-theoretic structure, and these are sets of quadratic forms. Dedekind developed his theory of ideals in order to supplement rings of algebraic integers with ``ideal divisors,'' extending the unique factorization property of the ordinary integers to these more general domains. Dedekind found that these ideal divisors could be identified with sets of elements in the original ring, now known as ``ideals.'' Like the characters, the ideals of a ring of algebraic integers bear an algebraic structure, and Dedekind was adamant that they should be treated as bona-fide objects.\footnote{See Avigad \cite{avigad:06}, especially page 172, and Edwards \cite{edwards:80}.} Similarly, Dedekind constructed the real numbers by identifying each of them with a pair of sets of rational numbers \cite{dedekind:72}. By the end of the century, it was common to view a quotient group as a group whose elements are equivalence classes, or cosets.\footnote{See the detailed discussion in Schlimm \cite[Section 3]{schlimm:08}. Other nice examples of pieces of nineteenth century mathematics that push in favor of set-theoretic abstraction are discussed in Wilson~\cite{wilson:10,wilson:unp}.}

The reasons given above to treat mathematical functions and sets as objects also speak in favor of treating them extensionally. The statements that ``there are two characters on $\ZZ / 2\ZZ$'' and ``there are $\ph(m)$ characters on $(\ZZ / m\ZZ)^*$'' are false if we take characters to be representations, as there are many different representations of the same character. We could, of course, develop notions of ``counting up to equivalence.'' In the early days of finite group theory, Camille Jordan described quotient groups as systems just like ordinary groups except that equality is replaced by an appropriate equivalence relation.\footnote{Again, see Schlimm \cite[Section 3]{schlimm:08}.} But, if we do that, mathematical statements become ``relativized'' to the appropriate equivalence relations, which constitute additional information that needs to be carried along and managed. The alternative is to extensionalize: then the only equivalence relation one has to worry about is equality.

We do not know the extent to which Frege was familiar with examples like these. But Wilson \cite{wilson:10,wilson:unp} calls attention to an important example of abstraction with which Frege was quite well acquainted. Frege was trained as a geometer, and studied under Ernst Schering in G\"ottingen. His dissertation, completed in 1873, was titled ``\"Uber eine geometrische Darstellung der imagin\"aren Gebilde in der Ebene'' (``On a geometric representation of imaginary forms in the plane''). Early nineteenth century geometers found great explanatory value and simplification in extending the usual Euclidean plane with various ideal objects, like ``points at infinity'' and ``imaginary'' points of intersection.  One of the few motivating examples that Frege provided in the \emph{Grundlagen} (\S64--\S68) is the fact that one can identify the ``direction'' of a line $a$ in the plane with the extension of the concept ``parallel to $a$.'' As Wilson points out (though Frege does not), these ``directions'' are exactly what is needed to serve as points at infinity, enabling one to embed the Euclidean plane in the larger projective plane, which has a number of pleasing properties. In the projective plane, all points have equal standing, and so it stands to reason that the concept-extensions used to introduce the new entities should be given the same ontological rights as the Euclidean points and lines used in their construction. Wilson characterizes such strategies for expansion as forms of ``relative logicism,'' since they provide a powerful means of relating the newly-minted objects to the more familiar ones.\footnote{See also Tappenden \cite{tappenden:06} for other ways that nineteenth century mathematics seems to have influenced his Frege's philosophical views.}

When it comes down to the nitty-gritty details, however, the only sustained formal development we have from Frege is his treatment of arithmetic. But even in this particular case, many of the issues we have raised come to the fore. In the \emph{Grundgesetze}, Frege defined a number of general operations and relations on tuples, sequences, functions, and relations. All of these now can be viewed as general set-theoretic constructions. What gives these constructions universal validity is that they can be applied to any domain of objects, and we now have great latitude in creating objects, as they are needed, to populate these domains. It is precisely the ability to bring a wide variety of mathematical constructions into the realm of objects, and the ability to define predicates and operations uniformly on this realm, that renders Frege's logic so powerful---too powerful, alas. But given Frege's goals, it should be clear why the extension operator held so much appeal.\footnote{The same uniformity is achieved in set theory by having a large universe of sets, and incorporating set-forming operations which return new elements of that universe. Russell introduced the notion of \emph{typical ambiguity} \cite{russell:08,feferman:82b} to allow ``polymorphic'' operations defined uniformly across types, and modern interactive theorem provers based on simple type theory follow such a strategy to obtain the necessary uniformities. For example, most such systems have operations $\mathit{card}_{\sigma}$ which maps a finite set of elements of type $\sigma$ to its cardinality, a natural number. The systems include mechanisms that allow one to define this family of operations uniformly, once and for all, treating $\sigma$ as a parameter. One can then write $\mathit{card} \; A$, and let the system infer the relevant type parameter from the type of $A$. This provides one means of coping with the nonuniformities that arise from a type-theoretic compartmentalization of the mathematical universe, but the difficulties that accrue to taking simple type theory as a mathematical foundation are complex; see, for example, \cite{avigad:12}.}

To sum up, we have traced a central tension in Frege's work to the need to balance two competing desiderata:
\begin{enumerate}
 \item the need for flexible but rigorous ways of talking about higher-type entities, like functions, predicates, and relations, without falling prey to incoherence; and
 \item the need for ways of dealing with mathematical objects uniformly, since mathematical constructions and operations have to be applied to many sorts of objects, many of which cannot be foreseen in advance.
\end{enumerate}
Compare this to the analysis of the mathematical pros and cons to treating functions as objects, as discussed in Section~\ref{analysis:section}, and note the similarities. 

At the end of an essay on Frege's treatment of concepts and objects, Thomas Ricketts briefly discusses aspects of nineteenth-century mathematical practice that may have had an influence on Frege. Agreeing with Wilson's assessment of the importance of being able to construct ideal elements in projective geometry, Ricketts writes:
\begin{quote}
Throughout his career, Frege is concerned with the introduction of new domains in mathematics, with the `creation' of new mathematical concepts. He vigorously polemicizes against formalist account of this practice and aims to develop an alternative to it. Frege's own approach here shines forth in a comment on Dedekind's account of the real numbers:
\begin{quote}
The most important thing for an arithmetician who recognizes in general the possibility of creation [of mathematical objects] will be to develop in an illuminating way [\emph{in einleuchender Weise}] the laws governing this in order to prove in advance of each individual creative act that the laws allow it. Otherwise, everything will be imprecise, and proofs will degenerate to a mere appearance, to a good-willed self-delusion. \cite[Vol.~2, \S 140]{frege:grundgesetze}.
\end{quote}
The desired foundation will be provided by formulating a logical law that, in the context of other logical laws, will yield as a theorem the existence of the desired new objects. \cite[pp.~217--218]{ricketts:10}
\end{quote}
What we have aimed to do here is to explain in greater detail why it is mathematically important to treat certain sorts of things as objects, and what, exactly, that amounts to. Not just Frege's work in geometry, but also his construction of the natural number system, would have impressed upon him the importance of having uniform operations and constructions on higher-order entities, and having a uniform way of making general assertions about these operations and constructions. 

In other words, Frege, like the various mathematical authors we have considered, was responding to methodological pressures that are inherent in the nature of the mathematical enterprise. As Ricketts emphasizes, Frege's entire foundational project was designed to address the important mathematical need of introducing clear means of expression, and developing general consensus as to the rules of use, while ensuring that the expressions and rules are meaningful, reliable, and consistent. While mathematicians from Dirichlet to Landau were focused on extending the edifice of mathematical knowledge, Frege's goal was to shore up the foundations. This difference translates to differences in perspective, focus, and method, but the distinctions are not sharp. Working from different ends of the spectrum, both Frege and his mathematical counterparts were working to clarify and extend mathematical method in powerful ways. In doing so, they addressed similar mathematical goals, and responded to similar mathematical constraints.

Frege is often faulted for failing to recognize the simple inconsistency that arises from the formal means he introduced to resolve the tension between the two concerns enumerated above. Nonetheless, it is worth highlighting the extent to which these two concerns were central to the subsequent development of logic and foundations. Russell's paradox shows that Frege was perfectly right to worry that an overly naive treatment of functions, concepts, and objects would lead to problems in the most fundamental use of our language and methods of reasoning. And, going into the twentieth century, developments in all branches of mathematics called for liberal means of constructing new mathematical domains and structures, as well as uniform ways of reasoning about their essential properties. The most fruitful and appropriate means of satisfying these needs was by no means clear at the turn of the twentieth century. Indeed, these issues were at the heart of the tumultuous foundational debates that were looming on the horizon.

\section*{Appendix: From cyclotomy to Dirichlet's theorem}
\addcontentsline{toc}{section}{Appendix: From cyclotomy to Dirichlet's theorem}
\label{appendix}

In Section~\ref{overview:section}, we sketched Dirichlet's approach to proving his theorem on primes in an arithmetic progression. Our goal here is to explain how Dirichlet is likely to have come upon his method of modifying Euler's argument to tease apart the contribution of the primes in each residue class from the overall sum of their reciprocals. 

Recall that if we split up the sum in Euler's equation (\ref{euler:primes:eqn}), we obtain
\begin{equation}
\log\sum_{n=1}^{\infty}\frac{1}{n^s}=\sum_{q \equiv 1 \bmod p} \frac{1}{q^{s}} + \sum_{q \equiv 2 \bmod p} \frac{1}{q^{s}} + \ldots + \sum_{q \equiv p-1 \bmod p} \frac{1}{q^{s}} + O(1).\tag{\ref{euler:primes:eqn:b}}
\end{equation}
As explained in Section~\ref{dirichlet:approach:section}, this shows that (\ref{euler:primes:eqn}) is too crude to prove Theorem~\ref{dirichlet:theorem}: we need to know that each of the terms on the right-hand side tends to infinity, not just their sum. It is here that ideas from the theory of equations are helpful. They come into play specifically in the theory of cyclotomy from Gauss' \emph{Disquisitiones Arithmeticae}, work with which Dirichlet was intimately acquainted. Historical overviews of the relevant ideas can be found in excellent books by Edwards and Tignol on the history of the theory of equations \cite{edwards:84,tignol:01}, and Curtis' equally impressive history of representation theory \cite{curtis:99}. Curtis also explains the role of characters in the theory of cyclotomy and Dirichlet's proof. What we aim to do here is make the progression of ideas leading from cyclotomy to Dirichlet's proof as explicit as possible. 

An important concern in the field of algebra is the extent to which the roots of a polynomial can be expressed in terms of arithmetic operations on the coefficients together with the extraction of roots. The quadratic formula dates to antiquity, and solutions to the cubic and quartic were presented by Cardano in his \emph{Ars Magna} of 1542. A natural challenge was then to determine a similar formula for the quintic. In 1770, Lagrange presented a general method of attacking this problem, using what has come to be known as the \emph{Lagrange resolvent}. Let $t_0, \ldots, t_{n-1}$ be the roots of the $n$th degree polynomial in question, and let $\omega$ be an $n$th root of unity, that is, a solution to the equation $\omega^n = 1$. Notice that 1 is always a solution to this equation, but there are $n - 1$ others. In fact, all of the roots can be taken to be powers a single ``primitive'' root of unity; for example, taking $\omega$ to be the complex number $e^{2 \pi i / n}$ will do. Lagrange considered the quantity
\[
 t_0 + \omega t_1 + \omega^2 t_2 + \ldots + \omega^{n-1} t_{n-1},
\]
as well as the quantities obtained by permuting the roots $t_0, \ldots, t_{n-1}$. Suppose $\omega$ is a primitive $n$th root of unity, and consider the values obtained by replacing $\omega$ in the previous expression with each of the values $1, \omega, \omega^2, \ldots, \omega^{n-1}$:
\begin{align*}
 x_0 & = t_0 + t_1 + t_2 + \ldots + t_n \\
 x_1 & = t_0 + \omega t_1 + \omega^2 t_2 + \ldots + \omega^{n-1} t_{n-1} \\
 x_2 & = t_0 + \omega^2 t_1 + \omega^4 t_2 + \ldots + \omega^{2(n-1)} t_{n-1} \\
 \vdots \\
 x_{n-1} & = t_0 + \omega^{n - 1} t_1 + \omega^{2 (n - 1)} t_2 + \ldots + \omega^{(n-1)^2} t_{n-1} 
\end{align*}
Lagrange observed that one can solve for each of $t_0, t_1, \ldots, t_{n-1}$ in terms of $x_0, x_1, \ldots, x_{n-1}$. For example, consider $x_0 + x_1 + \ldots + x_{n-1}$. Summing the first column gives $n \cdot t_0$. Summing the second column gives $t_1 \cdot (1 + \omega + \omega^2 + \ldots + \omega^{n-1})$. But because $\omega$ is a root of 
\[
 \omega^n - 1 = (\omega - 1) (\omega^{n-1} + \ldots + \omega^2 + \omega + 1)
\]
and $\omega \neq 1$, we have $1 + \omega + \omega^2 + \ldots + \omega^{n-1} = 0$. Similarly, summing the third column gives $t_2 \cdot (1 + \omega^2 + \omega^4 + \ldots + \omega^{2(n-1)})$; but $\omega^2$ is also an $n$th root of unity, and if $\omega$ is primitive (and $n > 2$), $\omega^2$ is also not equal to $1$, and the same argument shows that this quantity sums to $0$. The same argument shows that the remaining columns also sum to $0$, so we have $t_0 = (x_0 + \ldots + x_{n-1}) / n$, which is the desired expression for $t_0$. 

A similar trick works to compute the other values $t_k$: multiplying the $i$th equation by $\omega^{-ik}$ simply ``rotates'' the powers of $\omega$, leaving $1$'s in the $k$th column. Thus we have 
\[
t_k = \frac{1}{n} \sum_{i = 0}^{n-1} \omega^{-ik} x_i,
\]
which provides an expression for $t_k$ in terms of $t_0, \ldots, t_{n-1}$. Lagrange went on to consider the values of $x_0, \ldots, x_{n-1}$ that are obtained by replacing $\omega$ with other roots of unity, and conditions under which one can solve for those values, and hence $x_0, \ldots, x_n$, in terms of radicals. In doing so, he was analyzing and generalizing methods of solving equations developed by Vi\`ete, Tschrinhaus, and others who had come before. He showed that this ideas can be used to account for the known solutions to the quadratic, cubic, and quartic equations. 

The methods break down for the general solution to the quintic, but variations on the method can, however, be used to determine roots of \emph{particular} polynomials. Consider, for example, the polynomial $x^n - 1$ itself. We have already noted that we have $x^n - 1 = (x - 1) (x^{n-1} + x^{n-2} + \ldots + x^2 + x + 1)$. If $n$ is not a prime number, the second term can be factored into polynomials of lower degree, until one reaches polynomials that can no longer be factored; these are called \emph{irreducible} polynomials. The task of determining the roots of these polynomials is known as ``cyclotomy,'' or ``circle division,'' because the $n$ complex roots of $x^n - 1$ are evenly spaced around the unit circle in the complex plane.

The problem can be reduced to the case where $n$ is a prime number, which we will denote $p$ instead. In that case, $x^{p-1} + x^{p-2} + \ldots + x^2 + x + 1$ is irreducible, and if $\alpha$ is any root of this polynomial, the other $p-2$ roots are $\alpha^2, \alpha^3, \ldots, \alpha^{p-1}$. An expression for these roots in terms of radicals were provided by Vandermonde for the case where $p = 11$, and the general problem was taken up by Gauss in the last chapter of the \emph{Disquisitiones}.\footnote{Gauss was particularly interested in the case where $p$ is a prime number of the form $2^m - 1$, and showed that in that case, the solution enables one to carry out a geometric construction using compass and straightedge that divides the circle into $p$ equal parts. The \emph{Disquisitiones} hints as the solution to the general case, but both Edwards \cite{edwards:84} and Tignol \cite{tignol:01} observe that their are gaps in the presentation; a complete solution was provided by Galois.} The solution involves using the Lagrange resolvent, and taking the roots $t_0, t_1, \ldots, t_{p-2}$ to be the $p - 1$ roots $\alpha, \alpha^2, \ldots, \alpha^{p-1}$, but in a particular order.

The proof involves choosing, for the prime $p$ in question, a primitive element $g$ modulo $p$. Recall from Section~\ref{dirichlet:approach:section} that this means that the powers $g^0, g^1, \ldots, g^{p-1}$ modulo $p$ are exactly the nonzero residues modulo $p$. The solution to the equation $x^{p-1} + x^{p-2} + \ldots + x^2 + x + 1 = 0$ is obtained by considering the Lagrange resolvent
\begin{equation}
\label{cyclotomy:equation}
 \alpha^{g^0} \cdot \omega^0 + \alpha^{g^1} \cdot \omega^1 + \alpha^{g^2} \cdot \omega^2 + \ldots + \alpha^{g^{p-2}} \cdot \omega^{p-2},
\end{equation}
where $\omega$ is a $(p-1)$st root of unity. If we define $t_i$ to be $\alpha^{g^i}$, then this expression becomes
\[
 t_0 + \omega t_1 + \omega^2 t_2 + \ldots + \omega^{p-2} t_{p-2},
\]
and we are in the situation analyzed above. Lagrange's tricks tells us that if we can solve for the values of this expression when $\omega$ is replaced by $1, \omega, \omega^2, \ldots, \omega^{p-2}$ in succession, we can solve for all the values of $\alpha^{g^i}$, which are just the values $\alpha, \alpha^2, \ldots, \alpha^{p-1}$ written in a different order.

The reason for writing the elements $\alpha$ in the particular order they appear in (\ref{cyclotomy:equation}) is that, when they are written in that order, it \emph{is} possible to solve for each $t_i$, with an expression involving radicals. The details of the solution are not relevant to the proof of Dirichlet's theorem, but one particular aspect of the solution is. What makes the argument work is the careful pairing of $\alpha^{g^i}$ with $\omega^i$, which has the effect that for each $i$ and $j$, the element
\[
(\alpha^{g^i})^{g^j} = \alpha^{g^i \cdot g^j} = \alpha^{g^{i + j}}
\]
is paired with $\omega^{i + j}$. Using the notion of ``index'' defined in Section~\ref{dirichlet:approach:section}, we can express this as follows: for any $m$ and $n$, $\alpha^m$ is paired with $\omega^{\gamma_m}$, $\alpha^n$ is paired with $\omega^{\gamma_n}$, and $\alpha^{mn}$ is paired with $\omega^{\gamma_{mn}} = \omega^{\gamma_m + \gamma_n} = \omega^{\gamma_m} \omega^{\gamma_n}$. In other words, the key property used in the calculation of the roots of cyclotomic equations is that the map $m \mapsto \omega^{\gamma_m}$ is \emph{multiplicative} on the nonzero residues modulo $p$.

Dirichlet's great insight is that these ideas can be applied in the number-theoretic setting at hand, using the fact that the Euler product formula holds more generally with such a multiplicative function in the numerator. In the case where the common difference is a prime number $p$, if we choose a primitive root $g$ modulo $p$ and define $t_i = \sum_{q \equiv g^i \bmod p} 1/ q^s$, then equation (\ref{euler:character:eqn}) in Section~\ref{dirichlet:approach:section} can be written
\[
\log\sum_n\frac{\omega^{\gamma_n}}{n^s} = t_0 + t_1 \omega + \ldots + t_{p-2} \omega^{p-2} + O(1).
\]
The derivation of this equation relies on the generalized Euler formula, which requires that the map $m \mapsto \omega^{\gamma_m}$ is multiplicative. But once we have the equation in hand, we need only use the Lagrange trick, which is exactly what Dirichlet did. The more general case where $p$ is replaced by an arbitrary modulus $k$ is technically more difficult, but it builds on the same idea, combined with the behavior of the multiplicative group of residues modulo $k$ that are coprime to $k$. Once again, this is something which Dirichlet was intimately familiar with, from the work of Gauss.

\nocite{Fregecomp}


\begin{thebibliography}{99}

\bibitem{aspray:kitcher:88}
William Aspray and Philip Kitcher, editors. 
\newblock \emph{History and Philosophy of Modern Mathematics}.
\newblock University of Minnesota Press, Minneapolis, MN, 1988.

\bibitem{avigad:06}
Jeremy Avigad.
\newblock Methodology and metaphysics in the development of Dedekind's theory of ideals.
\newblock In \cite{ferreiros:gray:06}, pages 159--186. 

\bibitem{avigad:12}
Jeremy Avigad.
\newblock {Type inference in mathematics}.
\newblock {\em Bulletin of the European Association for Theoretical Computer
  Science}, 106:78--98, 2012.
  
\bibitem{avigad:unp}
Jeremy Avigad.
\newblock Modularity in mathematics.
\newblock {\em In preparation.}

\bibitem{avigad:morris:unp}
Jeremy Avigad and Rebecca Morris.
\newblock The concept of ``character'' in Dirichlet's theorem on primes in an arithmetic progression.
\newblock To appear in {\em Archive for History of Exact Sciences}. 

\bibitem{FregeReader}
Michael Beaney, editor.
\newblock {\em {The {F}rege Reader}}.
\newblock Blackwell Publishing, Malden, MA, 1997.

\bibitem{blanchette:12}
Patricia Blanchette.
\newblock Frege on shared belief and total functions.
\newblock {\em Journal of Philosophy},  CIX(1/2):9--39, 2012.

\bibitem{burgess:05}
John~P. Burgess.
\newblock {\em {Fixing {F}rege}}.
\newblock Princeton University Press,
  Princeton, NJ, 2005.

\bibitem{burgess:08}
John~P. Burgess.
\newblock {\em {Mathematics, Models, and Modality}}.
\newblock Cambridge University Press, Cambridge, 2008.

\bibitem{carnap:50}
Rudolf Carnap.
\newblock {Empiricism, semantics, and ontology}.
\newblock {\em Revue Internationale de Philosophie}, 4:20--40, 1950.
\newblock Reprinted in P.\ Benacerraf and H.\ Putnam, editors, \emph{Philosophy
  of Mathematics}, second edition, pages
  241--257, Cambridge University Press, Cambridge, 1983.

\bibitem{church:40}
Alonzo Church.
\newblock {A formulation of the simple theory of types}.
\newblock {\em Journal of Symbolic Logic}, 5:56--68, 1940.

\bibitem{curtis:99}
Charles W.~Curtis.
\newblock Pioneers of Representation Theory: Frobenius, Burnside, Schur, and Brauer.
\newblock American Mathematical Society and London Mathematical Society, 1999.

\bibitem{Poussin96}
Charles~Jean de~la Vall\'{e}e~Poussin.
\newblock D\'emonstration simplifi\'ee du th\'eor\`em de {D}irichlet sur la
  progression arithm\'etique.
\newblock {\em M\'{e}moires couronn\'{e}s et autres m\'{e}moires publi\'{e}s
  par L'Acad\'{e}mie Royale des Sciences, des Lettres et des Beaux-Arts de
  Belgique}, 53, 1895-1896.

\bibitem{dedekind:72}
Richard Dedekind.
\newblock {\em {Stetigkeit und irrationale Zahlen}}.
\newblock F.~Vieweg \& Sohn, Braunschweig, 1872.
\newblock Reprinted in \cite{dedekind:68}, volume 3, chapter L, pages 315--334.
  Translated by Wooster Beman as {``Continuity and irrational numbers''} in
  \emph{Essays on the Theory of Numbers}, Open Court, Chicago, 1901; reprinted
  by Dover, New York, 1963. The Beman translation is reprinted, with
  corrections by William Ewald, in \cite{ewald:96}, volume 2, pages 765-779.

\bibitem{dedekind:68}
Richard Dedekind.
\newblock {\em {Gesammelte mathematische {W}erke}}, volumes 1--3, edited by Robert Fricke, 
  Emmy Noether and \"{O}ystein Ore, F. Vieweg \& Sohn, Braunschweig, 1932. Reprinted by
  Chelsea Publishing Co., New York, 1968.

\bibitem{Dirichletshort}
Johann Peter Gustav Lejeune~Dirichlet.
\newblock {Beweis eines Satzes \"{u}ber die arithmetische Progression}.
\newblock {\em Bericht \"uber die Verhandlungen der k\"oniglich Presussischen Akademie der 
Wissenschaften Berlin}, 1837.
\newblock Reprinted in \cite{DirichletWerke1}, pages 309--312.
 
\bibitem{Dirichlet37}
Johann Peter Gustav Lejeune~Dirichlet Dirichlet.
\newblock {Beweis des Satzes, dass jede unbegrenzte arithmetische Progression,
  deren erstes Glied und Differenz ganze Zahlen ohne gemeinschaftlichen Factor
  sind, unendlich viele Primzahlen enth\"{a}lt}.
\newblock {\em Abhandlungen der k\"oniglich Preussischen Akademie der
  Wissenschaften}, pages 45--81, 1837.
\newblock Reprinted in \cite{DirichletWerke1}, pg 313--342. Translated by Ralf Stefan
  as ``There are infinitely many prime numbers in all
  arithmetic progressions with first term and difference coprime,''
  arxiv:0808.1408.
 
\bibitem{dirichlet:63b}
Johann Peter Gustave~Lejeune Dirichlet.
\newblock {\em {Vorlsesungen \"{u}ber {Z}ahlentheorie}}.
\newblock Vieweg, Braunschweig, 1863.
\newblock Edited by Richard Dedekind. Subsequent editions in 1871, 1879, 1894,
  with ``supplements'' by Richard Dedekind. Translated by
  John Stillwell, with introductory notes, as \emph{Lectures on Number Theory},
  American Mathematical Society, Providence, RI, 1999.

\bibitem{DirichletWerke1}
Johann Peter Gustav~Lejune Dirichlet.
\newblock {\em Werke},
\newblock edited by Leopold Kronecker. Georg Reimer, Berlin, 1889.

\bibitem{edwards:80}
Harold~M. Edwards.
\newblock {The genesis of ideal theory}.
\newblock {\em Archive for History of Exact Sciences}, 23:321--378, 1980.

\bibitem{edwards:84}
Harold~M. Edwards.
\newblock {\em Galois Theory}.
\newblock Springer, New York, 1984.

\bibitem{Euler48}
Leonhard Euler.
\newblock {\em {Introductio in analysin infinitorum, tomus primus}}.
\newblock Lausannae, 1748. Publications E101 and E102 in the Euler Archive.

\bibitem{ewald:96}
William Ewald, editor.
\newblock {\em {From Kant to Hilbert: A Source Book in the Foundations of
  Mathematics}}, volumes 1 and 2.
\newblock Oxford University Press, Oxford, 1996.
 
\bibitem{feferman:82b}
Solomon Feferman.
\newblock Typical ambiguity: trying to have your cake and eat it too.
\newblock In Godehard Link, editor, {\em One hundred years of {R}ussell's paradox}, pages 135--151, 
de Gruyter, Berlin, 2004.

\bibitem{ferreiros:gray:06}
Jos\'{e} Ferreir\'{o}s and Jeremy Gray, editors.
\newblock {\em {The Architecture of Modern Mathematics}}.
\newblock Oxford University Press, Oxford, 2006.

\bibitem{Frege84}
Gottlob Frege.
\newblock {\em {{D}ie {G}rundlage der {A}rithmetik}}.
\newblock W.~Koebner, Breslau, 1884.
\newblock Translated by J.~L.~Austin as \emph{The Foundations of Arithmetic}, 
  second revised edition, Northwestern University Press, Evanston, Illinois, 1980. 
  Excerpts translated by Michael Beaney in \cite{FregeReader}, pages 84--129.

\bibitem{Frege91}
Gottlob Frege.
\newblock {\em {F}unction und {B}egriff}.
\newblock Hermann Pohle, Jena, 1891.
\newblock Reprinted in \cite{Fregecomp}, and translated by Peter Geach as
  ``Function and concept'' in \cite{FregeReader}. Page number references are to the original publication.

\bibitem{cando}
Gottlob Frege.
\newblock {\"Uber Begriff und Gegenstand}.
\newblock {\em Vierteljahresschrift f\"{u}r wissenschaftliche Philosophie}, 16:192--505, 1892.
\newblock Reprinted in \cite{Fregecomp} and translated by Peter Geach as ``Concept and object'' 
in \cite{FregeReader}. Page number references are to the original publication.

\bibitem{frege:grundgesetze}
Gottlob Frege.
\newblock {\em Grundgesetze der arithmetik.}
\newblock H.~Pohle, Jena, volume 1, 1893, volume 2, 1903. Excerpts translated by Michael Beaney in \cite{FregeReader}, pages 194--223 and 258--289.

\bibitem{Frege04}
Gottlob Frege.
\newblock Was ist eine {F}unktion?
\newblock In Stefan Meyer, editor, {\em Festschrift Ludwig Boltzmann gewidmet
  zum sechzigsten Geburtstage}. J. A. Barth, Leipzig, 1904.
\newblock Reprinted in \cite{Fregecomp} and translated as ``What is a
  function?'' in Peter Geach and Max Black, \emph{Translations from the Philosophical Writings of Gottlob
  Frege}, Oxford University Press, Oxford, 1980. 
  Page number references are to the original publication.

\bibitem{Fregecomp}
Gottlob Frege.
\newblock {\em {F}unktion -- {B}egriff -- {B}edeutung}, edited by Mark Textor.
\newblock Vandenhoeck and Ruprecht, G\"ottingen, 2002.

\bibitem{gauss:01}
Carl~Friedrich Gauss.
\newblock {\em {Disquisitiones Arithmeticae}}.
\newblock G. Fleischer, Leipzig, 1801.
\newblock Reprinted in Gauss' \emph{Werke}, K\"{o}niglichen Gesellschaft der
  Wissenschaften, G\"{o}ttingen, 1863. Translated with a preface by Arthur A.
  Clarke, Yale University Press, New Haven, 1966, and republished by Springer, 
  New York, 1986.

\bibitem{gordon:melham:93}
M.~J.~C. Gordon and T.~F. Melham, editors.
\newblock {\em {Inroduction to HOL: A Theorem Proving Environment for
  Higher-Order Logic}}.
\newblock Cambridge University Press, Cambridge, 1993.

\bibitem{gray:92}
Jeremy Gray.
\newblock The nineteenth-century revolution in mathematical ontology.
\newblock
In Donald Gillies, editor, \emph{Revolutions in Mathematics}, 226--248, Oxford University Press, 
Oxford, 1992.

\bibitem{harrison:07c}
John Harrison.
\newblock {{HOL} light: a tutorial introduction}.
\newblock In Mandayam Srivas and Albert Camilleri, editors, {\em {Proceedings
  of the First International Conference on Formal Methods in Computer-Aided
  Design}}, pages 265--269, Springer, Berlin, 1996.

\bibitem{hawkins:71}
Thomas Hawkins.
\newblock {The origins of the theory of group characters}.
\newblock {\em Archive for History of Exact Sciences}, 7:142--170, 1971.

\bibitem{heck:97}
Richard Heck.
\newblock The {J}ulius {C}easar objection.
\newblock In Richard Heck, editor, \emph{Language, Thought, and Logic: 
Essays in Honour of Michael Dummett}, pages 273--308, Oxford University Press, Oxford, 1997.
\newblock Reprinted in Richard Heck, editor, \emph{Frege's Theorem}, pages 127--155, Oxford
  University Press, Oxford, 2011.

\bibitem{compthought}
Richard Heck and Robert May.
\newblock The composition of thoughts.
\newblock {\em No\^{u}s}, 45:126--166, 2011.

\bibitem{kitcher:88}
Philip Kitcher.
\newblock {Mathematical naturalism}.
\newblock In \cite{aspray:kitcher:88}, pages 293--325.

\bibitem{kronecker:70}
Leopold Kronecker.
\newblock Auseinandersetzung einiger eigenschaften der klassenzahl idealer
  complexer zahlen.
\newblock {\em Monatsberichte der K\"oniglich Preussischen Akademie der
  Wissenschaften zu Berlin}, pages 881--882, 1870.
\newblock Reproduced in \cite{kronecker:68}, volume I, pages 271--282.

\bibitem{kronecker:68}
Leopold Kronecker.
\newblock {\em {W}erke}, edited by Kurt Hensel, volumes 1--5,
B.~G.~Teubner, Leipzig, 1895--1930. 
\newblock Reprinted by Chelsea Publishing Co., New York,
  1968.

\bibitem{kummer:46}
Ernst~Eduard Kummer.
\newblock {Zur Theorie der complexen Zahlen}.
\newblock {\em Koniglich Akademie der Wissenschaft Berlin, Monatsbericht},
  pages 87--97, 1846.
\newblock Also in \emph{Journal f\"ur die reine und angewandte
  Mathematik} 35:319--326, 1847, and in Kummer's \emph{Collected
  Papers}, edited by Andr\'e Weil, Springer-Verlag, Berlin, 1975, volume 1,
  203--210.
  
\bibitem{landau1909}
Edmund Landau.
\newblock {\em {Handbuch der {L}ehre von der {V}erteilung der {P}rimzahlen}},
  volume~1.
\newblock B.~G.~Teubner, Leipzig, 1909.
 
\bibitem{mach:93}
Ernst Mach.
\newblock {\em {The Science of Mechanics: {A} Critical and Historical Account
  of its Development}}, translated by Thomas J. McCormack. The Open Court Publishing
  Co., La Salle, Illinois, 1960.

\bibitem{maddy:97}
Penelope Maddy.
\newblock {\em {Naturalism in Mathematics}}.
\newblock Oxford University Press, New York, 1997.

\bibitem{mancosu:08}
Paolo Mancosu, editor. 
\newblock {\em {The Philosophy of Mathematical Practice}}.
\newblock Oxford University Press, Oxford, 2008.
 
\bibitem{manders:08}
Kenneth Manders.
\newblock {The {E}uclidean diagram}.
\newblock In \cite{mancosu:08}, pages 80--133. 

\bibitem{morris:11}
Rebecca Morris.
\newblock \emph{Character and object}.
\newblock Master's thesis, Carnegie Mellon University, 2011.

\bibitem{nipkow:et:al:02}
Tobias Nipkow, Lawrence~C. Paulson, and Markus Wenzel.
\newblock {\em {Isabelle/{HOL}. A Proof Assistant for Higher-Order Logic}}.
\newblock Springer, Berlin, 2002.

\bibitem{panza:unp}
Marco Panza.
\newblock {From {L}agrange to {F}rege: functions and expressions}.
\newblock In Carlos Alvarez and Andrew Arana, editors, 
\emph{Analytic Philosophy and the Foundations of Mathematics}, 
Palgrave Macmillan, to appear. 

\bibitem{parsons:08}
Charles Parsons.
\newblock \emph{Mathematical Thought and its Objects}.
\newblock Cambridge University Press, Cambridge, 2008.

\bibitem{putnam:12}
Hilary Putnam.
\newblock \emph{Philosophy in an Age of Science: Physics, Mathematics, and Skepticism}.
\newblock Harvard University Press, Cambridge, 2012.

\bibitem{quine:48}
W.~V.~O. Quine.
\newblock On what there is.
\newblock {\em The Review of Metaphysics}, 2:21--38, 1948.
\newblock Reprinted in W.~V.~O.~Quine, \emph{From a Logical Point of View},
  Harvard University Press, Cambridge, 1980.
% 
\bibitem{quine:51}
W.~V.~O. Quine.
\newblock {Main trends in recent philosophy: two dogmas of empiricism}.
\newblock {\em The Philosophical Review}, 60:20--43, 1951.

\bibitem{quine:55}
W.~V.~O. Quine.
\newblock {On {F}rege's way out}.
\newblock {\em Mind}, 64:145--159, 1955.

\bibitem{quine:69}
W.~V.~O. Quine.
\newblock {\em {Ontological Relativity, and other Essays}}.
\newblock Columbia University Press, New York, 1969.
% 
\bibitem{quine:60}
W.~V.~O. Quine.
\newblock {\em {Word and Object}}.
\newblock MIT Press, Cambridge, MA, 1960.
  
\bibitem{reck:07}
Erich Reck. 
\newblock Frege-Russell numbers: analysis or explication?
\newblock In Michael Beaney, editor, \emph{The Analytic Turn: Analysis in 
Early Analytic Philosophy and Phenomenology}, pages 33-50, 
Routledge, London, 2007.

\bibitem{ricketts:10}
Thomas Ricketts. 
\newblock Concepts, objects, and the context principle.
\newblock In Michael Potter and Tom Ricketts, eds., {\em The Cambridge Companion to Frege}.
Cambridge University Press, Cambridge, 2010, pages 149--219.
 
\bibitem{russell:08}
Bertrand Russell.
\newblock Mathematical logic as based on the theory of types.
\newblock {\em American Journal of Mathematics}, 30:222--262, 1908.

\bibitem{russell:whitehead:10}
Bertrand Russell and Alfred~North Whitehead.
\newblock {\em {Principia Mathematica}}.
\newblock Cambridge University Press, Cambridge, first volume, 1910; second volume, 1912;
  third volume, 1913.

\bibitem{schlimm:08}
Dirk Schlimm.
\newblock {On abstraction and the importance of asking the right research
  questions: could {J}ordan have proved the {J}ordan-{H}\"{o}lder theorem?}
\newblock {\em Erkenntnis}, 68:409--420, 2008.
% 
% \bibitem{stein:88}
% Howard Stein.
% \newblock {Logos, logic, and logistik\'{e}}.
% \newblock In \cite{aspray:kitcher:88}, pages 238--259.

\bibitem{tait:86}
William Tait.
\newblock Truth and proof: the Platonism of mathematics.
\newblock \emph{Synthese}, 69: 341--370. Reprinted in William Tait, \emph{The Provenance of Reason: Essays in The Philosophy of Mathematics and its History}, Oxford University Press, Oxford, 2005, pages 61--88.

\bibitem{tappenden:06}
Jamie Tappenden.
\newblock {The {R}iemannian background to {F}rege's philosophy}.
\newblock In \cite{ferreiros:gray:06}, pages 97--132.

\bibitem{tignol:01}
Jean-Pierre Tignol.
\newblock {\em Galois' Theory of Algebraic Equations.}
\newblock World Scientific, New Jersey, 2001.

\bibitem{urquhart:08}
Alasdair Urquhart.
\newblock Mathematics and physics: strategies of assimilation.
\newblock In \cite{mancosu:08}, pages 417--440.

\bibitem{wilson:94}
Mark Wilson.
\newblock Can we trust logical form?
\newblock {\em The Journal of Philosophy}, 91:519--544, 1994.

\bibitem{wilson:10}
Mark Wilson.
\newblock Frege's mathematical setting.
\newblock In Michael Potter and Tom Ricketts, editors, {\em The Cambridge Companion to Frege}, 
Cambridge University Press, Cambridge, 2010, page 379--412.

\bibitem{wilson:unp}
Mark Wilson.
\newblock Enlarging one's stall \emph{or} How did all of these \emph{sets} get in here?
\newblock Preprint.
 
\bibitem{wittgenstein:89}
Ludwig Wittgenstein.
\newblock {\em Wittgenstein's Lectures on the Foundations of Mathematics,
  Cambridge, 1939}.
\newblock University of Chicago Press, Chicago, 1989.

\end{thebibliography}
\end{document}